# Efficient Dynamic Allocation Policy for Robust Ranking and Selection under Stochastic Control Framework

Hui Xiao[1], Zhihong Wei[1]

[1]School of Management Science and Engineering, Southwestern University of Finance and Economics, 555 Liutai Avenue, Wenjiang District, Chengdu 611130, P. R. China.

**Abstract:** This research considers the ranking and selection with input uncertainty. The objective is to maximize the posterior probability of correctly selecting the best alternative under a fixed simulation budget, where each alternative is measured by its worst-case performance. We formulate the dynamic simulation budget allocation decision problem as a stochastic control problem under a Bayesian framework. Following the approximate dynamic programming theory, we derive a one-step-ahead dynamic optimal budget allocation policy and prove that this policy achieves consistency and asymptotic optimality. Numerical experiments demonstrate that the proposed procedure can significantly improve performance.

**Keywords:** Simulation optimization, Ranking and selection, Stochastic control problem, Bayesian, Input uncertainty.

## 1. Introduction

Decision-making processes are widespread in many discrete-event dynamic systems (DEDS), such as DEDS in manufacturing, transportation, and healthcare management. However, as the complexity of DEDS gradually increases, the assumptions in traditional analytical models are hardly satisfied. Thus, modeling DEDS using analytical models is generally infeasible, which poses a significant challenge for evaluating the decision-making of DEDS. (Cassandras and Lafortune 2009). Simulation is a powerful modeling tool that uses logically tricky, often non-mathematical models and few assumptions to analyze DEDS and evaluate system decisions effectively. (Fu et al., 2015). Nevertheless, system simulation is usually expensive and time-consuming, especially when the system involved is complex (Ho et al., 2007). Meanwhile, the stochastic nature of simulation models leads decision-makers to require large simulation budgets with multiple simulation replications for each alternative to obtain stable system performance estimates (Chen & Lee, 2011). Therefore, how to optimize the simulation budget allocation to improve the efficiency of system simulation is one of the current popular topics of simulation optimization.

The ranking and selection (R&S) procedures are classical methods designed to improve the efficiency of simulations, which intelligently allocate the simulation budget to each alternative based on information, such as the sample mean and sample variance, to select the best alternative from among all alternatives efficiently. The R&S procedures can be divided into frequentist and Bayesian streams. (Kim & Nelson, 2006; Hunter & Nelson, 2017; Hong et al., 2021). The first stream focused on finding a feasible simulation budget allocation rule to guarantee that the pre-specified probability of correct selection (PCS) is achieved, such as the indifference zone (IZ, e.g., Rinott 1978; Frazier 2014, and Fan 2016). Whereas the latter aims to achieve the PCS as high as possible under the same simulation budget, such as optimal computing budget allocation (OCBA, Chen et al. 2000), knowledge gradient (KG, Frazier, et al. 2009), expected value information (EVI, Chick et al. 2010), expected improvement (EI, Ryzhov 2016), and asymptotically optimal allocation procedure (AOAP, Peng et al. 2018).

So far, the R&S procedures have been further extended and used to solve many simulation optimization problems. Such as the subset selection problem (e.g., Chen et al., 2008; Zhang et al., 2015; Gao & Chen, 2016), the multi-objective simulation optimization problem (Lee et al., 2010), the stochastic simulation time optimization problem (Jia, 2012), the best alternative selection problem under different contexts (Cakmak, 2021; Li et al., 2022), the best quantile selection problem (Peng et al., 2021), the covariate simulation optimization problem (Shen et al., 2021), the data censored simulation optimization problem (Liu et al., 2023), and the large scale simulation optimization problem (Zhong et al., 2022). In addition, the application of R&S procedures in simulation can be found in electric vehicle charging management (Jiang et al., 2020), transportation (Xiao et al., 2021), healthcare management (Fan et al., 2020; Li et al., 2022), etc.

In this work, we consider how to optimize the allocation of a fixed simulation budget $T$ to maximize the probability of correctly selecting the best alternative with input uncertainty

($PCS_W$), where their performance is measured by the worst-case performance over a finite set of possible distributions. The most relevant work to this paper is Gao et al. (2017), which constructed the above selection of the best problem with input uncertainty under a fixed simulation budget constraint based on the OCBA method, derived an optimal allocation rule for the simulation budget of this selection problem, and proposed an optimal computing budget allocation procedure (ROCBA) according to the allocation rule. The allocation decision for the simulation budget is a stochastic control problem. However, the allocation rules used by ROCBA are derived using static optimization methods under the assumption of perfect information (e.g., assuming that the mean and variance are known). Therefore, in some low-confidence scenarios (i.e., the means of all alternatives are close, the variances are relatively large, and the simulation budgets are relatively limited), the static allocation rules used by ROCBA may inadequately determine the dynamic allocation decision behavior of the simulated budget, thus making the algorithm appear to perform poorly (Chen et al., 2006; Peng et al., 2017;). Conversely, it would be beneficial to model and analyze the problem using stochastic dynamic programming methods that accurately describe the dynamic simulation budget allocation process (Peng et al., 2016; Peng et al., 2018).

Based on the above analysis, this paper addresses the problem of selecting the best alternative with input uncertainty from the perspective of stochastic dynamic programming. Firstly, we formulate the dynamic simulation budget allocation decision process as a stochastic control problem in the Bayesian framework. Then, with the approximation of $PCS_W$ as the value function, and the approximation is approximated one step ahead according to the approximate dynamic programming approach. Thus, we propose an efficient asymptotically optimal dynamic simulation budget allocation policy and its corresponding one-step-ahead algorithm for selecting the best alternative. Finally, we prove that the allocation policy can achieve consistency. That is, the best alternative selected by the allocation policy is true when the simulation budget goes to infinity, and the simulation allocation ratio generated by the allocation policy converges asymptotically to the optimal simulation budget allocation ratio proposed by Gao et al., (2017) based on OCBA Method. Numerical experiments demonstrate the efficiency of the proposed procedure.

The rest of the paper is organized as follows. Section 2 describes the selection of the best problem with input uncertainty, formulates the dynamic sampling decision problem as the stochastic control problem, and gives assumptions for this paper. Section 3 derives the dynamic sampling policy. Section 4 analyzes the properties of the dynamic sampling policy. Numerical experiments are provided in section 5 to test the performance of the proposed algorithm. The last section summarizes this paper.

## 2. Problem Formulation

There are $k$ competing alternatives and the set of these alternatives is $\mathcal{L} = \{\ell_1, \ell_2, \ldots, \ell_k\}$. Let $E_P[h(\ell, \xi)]$ be the performance measure function of alternative $\ell$ and must be estimated via Monte Carlo simulation. $h(\ell, \xi)$ is a random estimate of the performance of alternative $\ell \in \mathcal{L}$, where $\xi$ is the uncertainty parameter and follows an unknown distribution (scenario) $P$. In real-world applications, the parameters of scenario $P$ and the performance value $h(\ell, \xi)$ are estimated from historical data, which leads to the uncertainty of $P$. In order to characterize this input uncertainty, this work follows the literature with input uncertainty (Fan et al., 2013; Fan et al., 2020; Gao et al., 2017) and assumes that for $\forall \ell_i \in \mathcal{L}$, the scenario set $\mathcal{P}$ contains a finite number of $m$ possible scenarios of parameter $\xi$, i.e., $\mathcal{P} = \{P_1, \ldots, P_m\}$, where $\mathcal{P}$ incorporate the uncertainty from both the input distribution and its associate parameters.

To facilitate the following problem formulation, we first introduce some notations and assumptions. Let $i$ reprensents alternative $\ell_i$, $l$ is scenario $P_l$ and $(i, d)$ is alternative-scenario pair $(\ell_i, P_d)$ where $i \in \{1, 2, \ldots, k\}$ and $d \in \{1, 2, \ldots, m\}$. Further, let $X_{i,d,t}, t \leq T, t \in \mathbb{Z}^+$ is a simulation output of $h(\ell_i, \xi)$ with $\xi$ following distribution $P_d$ in $t$-th simulation replication, $\mu_{i,d} = E[X_{i,d,t}]$, and $\sigma_{i,d}^2 = \text{Var}[X_{i,d,t}]$, where $T < \infty$ is total number of simulation budget. We measure alternative $i$ by its worst-case mean performance $\mu_{i,d_i}$ over $\{1, 2, \ldots, k\}$, i.e., $\mu_{i,d_i} = \min_{l \in \{1,2,\ldots,m\}} \mu_{i,d}, i \in \{1, 2, \ldots, k\}$ where $d_i$ is the worst-case scenario of alternative $i$. Moreover, we further make the following assumption are need throughout the paper.

**Assumption 1.** The simulation output $X_{i,d,t}, i \in \{1, 2, \ldots, k\}, d \in \{1, 2, \ldots, m\}, t \leq T$ be independent and identical distribution (IID) from replication to replication and for alternative $i$ in scenario $j$.

This assumption is a canonical condition in R&S (Chen & Lee, 2011). Common random numbers and correlated sampling were not considered in this work.

**Assumption 2.** There always exists the worst-case scenario $d_i$ for each alternative $i \in \{1,2,\ldots,k\}$ such that $\mu_{i,d_i} < \mu_{i,d}, d_i \in \{1,2,\ldots,m\}, d_i \neq d$ and there exists alternative $i \in \{1,2,\ldots,k\}$ such that $\mu_{i,d_i} > \mu_{j,d_j}, j \in \{1,2,\ldots,m\}, i \neq j$.

Assumption 2 ensures that the worst-case scenario exists for each alternative and that the true best alternative is unique (Gao et al., 2017; Xiao et al., 2020).

Our aim is to develop an asymptotic optimal simulation budget allocation policy to correctly select the best alternative $\langle 1 \rangle$ from $\{1,2,\ldots,k\}$ with the largest worst-case mean performance under fixed simulation budget $T$, i.e.,

$$\langle 1 \rangle = \arg\max_{i \in \{1,2,\ldots,k\}} \mu_{i,d_i}. \tag{1}$$

Let $Q(\cdot;\theta)$ to be the joint distribution of $X_t = (X_{1,1,t}, \ldots, X_{1,m,t}, \ldots, X_{k,1,t} \ldots, X_{k,m,t})$, i.e., $X_t \sim Q(\cdot;\theta)$, where $\theta$ is a vector comprising of all unknown parameter in parametric space $\Theta$, and especially, $(\mu_{1,1}, \ldots, \mu_{1,m}, \ldots, \mu_{k,1}, \ldots, \mu_{k,m}) \in \theta$. Suppose $X_{i,d,t}$ follows a marginal distribution $Q_{i,d}(\cdot;\theta_{i,d})$, $i \in \{1,2,\ldots,k\}$ and $d \in \{1,2,\ldots,m\}$ where $\theta_{i,d}$ comprising of all unknown parameter in marginal distribution. Then, from Assumption 2, we have $Q(\cdot;\theta) = \prod_{i=1}^{k} \prod_{d=1}^{m} Q_{i,d}(\cdot;\theta_{i,d})$. Moreover, suppose the unknown parameter $\theta$ follows a prior distribution $F(\cdot;\varepsilon_0)$, i.e., $\theta \sim F(\cdot;\varepsilon_0)$, where prior information $\varepsilon_0$ comprises all hyper-parameter for the parametric space of prior distribution.

The sequence mappings $A_t(\cdot) = (a_1(\cdot),\ldots,a_t(\cdot))$ is simulation budget allocation policy where $a_t(\mathcal{F}_{t-1}^a) = (i,d) \in \{1,2,\ldots,k\} \times \{1,2,\ldots,m\}$, $0 < t < T$, which means that allocate the $t$-th simulation replication to alternative $(i,d)$ based on information set $\mathcal{F}_{t-1}^a$. The information set form as follows:

$$\mathcal{F}_t^a \triangleq \{a_t(\mathcal{F}_{t-1}^a); \mathcal{F}_t\}, \tag{2}$$

where $\mathcal{F}_t$ contains prior information $\varepsilon_0$ and all previous sample information. Let $t_{i,d} \triangleq \sum_{l=1}^{t} a_{i,d,l}(\mathcal{F}_{l-1}^a)$ as the simulation budget obtained for $(i,d) \in \{1,2,\ldots,k\} \times \{1,2,\ldots,m\}$ after $t$ simulated budget allocation, where $a_{i,d,l}(\mathcal{F}_{l-1}^a) \triangleq \mathbb{I}(a_{i,d}(\mathcal{F}_{l-1}^a) = i)$, and $\mathbb{I}(\cdot)$ is an indicator function that equal 1 if the event in the bracket. $X_{i,d}^{(t)} \triangleq (\tilde{X}_{i,d,1}, \ldots, \tilde{X}_{i,d,t_{i,d}})$ is the sample information of $(i,d)$ obtained by $t_{i,d}$ simulation replications, where $\tilde{X}_{i,d,\tau} = X_{i,d,l}$ when

$a_{i,d,l}(\mathcal{F}^a_{l-1}) = 1$, $\tau = \{1,2,\ldots,t_{i,d}\}$ and $l = \{1,2,\ldots,t\}$. Then, we can get $\mathcal{F}_t = \{\varepsilon_0, X^{(t)}_{1,1}, \ldots, X^{(t)}_{1,m}, \ldots, X^{(t)}_{k,1}, \ldots, X^{(t)}_{k,m}\}$, where $0 < t \leq T$ and $\mathcal{F}_0 = \varepsilon_0$.

Following the Bayesian rule, the posterior distribution of the unknown parameter $\theta$ conditioned on the information set $\mathcal{F}^a_t$ can take the following forms:

$$F(d\theta|\mathcal{F}^a_t) = \frac{L(\mathcal{F}^a_t;\theta)F(d\theta;\varepsilon_0)}{\int_{\theta\in\Theta}L(\mathcal{F}^a_t;\theta)F(d\theta;\varepsilon_0)}, \quad (3)$$

and the predictive distribution of $X_{i,d,t}$

$$Q_{i,d}(dx_{i,d,t+1}|\mathcal{F}^a_t) = \frac{\int_{\theta\in\Theta}Q_{i,d}(dx_{i,d,t+1};\theta_i)L(\mathcal{F}^a_t;\theta)F(d\theta;\varepsilon_0)}{\int_{\theta\in\Theta}L(\mathcal{F}^a_t;\theta)F(d\theta;\varepsilon_0)}, \quad (4)$$

where $(i,d) \in \{1,2,\ldots,k\} \times \{1,2,\ldots,m\}$, $0 < t < T$, $L(\cdot)$ is the likelihood function of $X_t = (X_{1,1,t},\ldots,X_{1,m,t},\ldots,X_{k,1,t}\ldots,X_{k,m,t})$, and $d(\cdot)$ is referred to as the Lebesgue (counting) measure for continuous (discrete) distribution.

Define the posterior best alternative $\langle 1 \rangle_t$ after $t$ simulation replications, where $\langle i \rangle_t \in \{1,2,\ldots,k\}$ is indices by largest posterior worst-case mean i.e., $\mu^{(t)}_{\langle 1 \rangle_t,d_{\langle 1 \rangle_t}} > \mu^{(t)}_{\langle 2 \rangle_t,d_{\langle 2 \rangle_t}} > \cdots > \mu^{(t)}_{\langle k \rangle_t,d_{\langle k \rangle_t}}$, where $d_{\langle i \rangle_t} = \underset{l\in\{1,2,\ldots,m\}}{\arg\min}\,\mu^{(t)}_{i,d}$ is the posterior worst-case scenario of alternative, $0 < t \leq T$. Moreover, we define the fixed selection policy $S(\mathcal{F}^a_T) = \langle 1 \rangle_T$ which makes the final selection after $T$-th simulation replication. A correct selection policy occurs when $\langle 1 \rangle_T = \langle 1 \rangle$ and the reward of selection policy can be measure by

$$PCS_W = P\{\mathbb{I}\{\langle 1 \rangle_T = \langle 1 \rangle\}|\mathcal{F}^a_T\} = P\left\{\bigcap_{i=1,i\neq 1}^{k}\bigcup_{d=1}^{m}\bigcap_{l=1}^{m}(\mu_{\langle 1 \rangle_T,l} \geq \mu_{\langle i \rangle_T,d})\middle|\mathcal{F}^a_T\right\}. \quad (5)$$

Generally, the posterior $PCS_W$ (5) has no closed form and can be estimated by Monte Carlo simulation. However, Monte Carlo simulation is expensive and time-consuming. Since simulation budget allocation aims to improve simulation efficiency, we need a relatively fast and inexpensive method to estimate (5). We derive the low bound of (5) as follows:

$$P\left\{\bigcap_{i=1,i\neq 1}^{k}\bigcap_{l=1}^{m}(\mu_{\langle 1 \rangle_T,l} \geq \mu_{\langle i \rangle_T,d_{\langle i \rangle T}})\middle|\mathcal{F}^a_T\right\}. \quad (6)$$

The (6) shows that (5) tends to be one when (6) tends to be one. Consequently, we use (6) to estimate (6) and maximize (5) by maximizing (6). So, we focus primarily on developing an asymptotically optimal dynamic simulation budget allocation policy $A$, where $A \triangleq A_T$, for maximizing (6) in the remaining part of this paper.

Under the Bayesian framework, we formulate the dynamic decision process in robust ranking and selection as a stochastic control problem (SCP) as follows:

The expected pay-off for the allocation policy $A$ can be recursively given by

$$V_T\left(\mathcal{F}_T^a; A\right) \triangleq \mathrm{P}\left\{\bigcap_{i=1, i\neq 1}^{k} \bigcap_{l=1}^{m}\left(\mu_{\langle 1\rangle_T, l} \geq \mu_{\langle i\rangle_T, d_{\langle i\rangle_T}}\right) \middle| \mathcal{F}_T^a\right\}, \tag{7}$$

and $0 \leq t < T$,

$$\begin{aligned}V_t\left(\mathcal{F}_t^a; A\right) &\triangleq \mathrm{E}\left[V_{t+1}\left(\mathcal{F}_t^a \cup \{X_{i,d,t+1}\}; A\right) \middle| \mathcal{F}_t^a\right]\bigg|_{(i,d)=a_{t+1}(\mathcal{F}_t^a)} \\ &= \int_{\chi_i} V_{t+1}\left(\mathcal{F}_t^a \cup \{x_{i,d,t+1}\}; A\right) Q_{(i,d)}\left(dx_{i,d,t+1} \middle| \mathcal{F}_t^a\right)\bigg|_{(i,d)=a_{t+1}(\mathcal{F}_t^a)},\end{aligned} \tag{8}$$

where $\chi_i$ is the support of $X_{i,d,t+1}$. Then, the optimal allocation policy $A^*$ as

$$A^* \triangleq \sup_{A} V_0(\varepsilon_0; A). \tag{9}$$

**Theorem 1.** Under the Bayesian framework, the information set $\mathcal{F}_T$ wholly determine the posterior distribute (3) of $\theta$, and the information set $\mathcal{F}_t$ wholly determine the predictive distribute (4) of $X_{i,d,t}$, and they are independent of the allocation policy $A$.

*Proof:* Since only $a_t(\mathcal{F}_{t-1}^a) = (i,d)$ obtain one simulation replication in any $t$ step, while the simulation replication of other alternative-scenario pairs are missing. Therefore, the likelihood of observations collected by the sequential simulation replication process though $t$ steps is given by

$$\begin{aligned}L\left(\mathcal{F}_t^a; \theta\right) &\triangleq \int \cdots \int_{\chi^t} \prod_{\tau=1}^{t} q(x_\tau; \theta) \prod_{i=1}^{k} \prod_{d=1}^{m} \left\{a_{i,d,\tau}\left(\mathcal{F}_{\tau-1}^a\right) \delta_{X_{i,d,\tau}}\left(dx_{i,d,\tau}\right)\right. \\ &\quad \left. + \left(1 - a_{i,d,\tau}\left(\mathcal{F}_{\tau-1}^a\right)\right) dx_{i,d,\tau}\right\} \\ &= \left[\sum_{i=1}^{k}\sum_{d=1}^{m} a_{i,d,\tau}\left(\mathcal{F}_{\tau-1}^a\right) q_{i,d}\left(X_{i,d,\tau}; \theta_{i,d}\right)\right] \int \cdots \int_{\chi^{t-1}} \prod_{\tau=1}^{t-1} q(x_{\tau-1}; \theta) \\ &\quad \times \prod_{i=1}^{k}\prod_{d=1}^{m}\left\{a_{i,l,\tau}\left(\mathcal{F}_{\tau-1}^a\right)\delta_{X_{i,d}}\left(dx_{i,d,\tau}\right) + \left(1 - A_{i,\tau}\left(\mathcal{F}_{\tau-1}^a\right)\right)dx_{i,a,\tau}\right\} \\ &= \prod_{\tau=1}^{t}\left[\sum_{i=1}^{k}\sum_{d=1}^{m}(a_{i,d,\tau}\left(\mathcal{F}_{\tau-1}^a\right)q_{i,l}\left(X_{i,d,\tau}; \theta_{i,d}\right))\right] = \prod_{i=1}^{k}\prod_{d=1}^{m}\prod_{\tau=1}^{t_{i,d}} q_{i,d}\left(\tilde{X}_{i,d,\tau}; \theta_{i,d}\right).\end{aligned} \tag{10}$$

in which $\chi \triangleq \chi_{1,1} \times ... \times \chi_{1,m} \times ... \times \chi_{k,1} \times ... \times \chi_{k,m}$, and $\delta_{i,d}(\cdot)$ is the delta measure with a mass point $x$. Observing (9), we can get that the information set $\mathcal{F}_t$ wholly determine $L(\mathcal{F}_t^a|\theta)$, so $L(\mathcal{F}_t|\theta) = L(\mathcal{F}_t^a|\theta)$. From (3)-(4), we have the posterior distribution of $\theta$ is

$$F\left(d\theta\mid \mathcal{F}_T^a\right)=\frac{L\left(\mathcal{F}_T;\theta\right)F\left(d\theta;\varepsilon_0\right)}{\int_{\theta\in\Theta}L\left(\mathcal{F}_T;\theta\right)F\left(d\theta;\varepsilon_0\right)}$$

$$=\frac{\prod_{i=1}^{k}\prod_{d=1}^{m}\prod_{\tau=1}^{t_{i,d}}qd\left(\tilde{X}_{i,d,\tau};\theta_{i,d}\right)F\left(d\theta;\varepsilon_0\right)}{\int_{\theta\in\Theta}\prod_{i=1}^{k}\prod_{d=1}^{m}\prod_{\tau=1}^{t_{i,d}}q_{i,d}\left(\tilde{X}_{i,d,\tau};\theta_{i,d}\right)F\left(d\theta;\varepsilon_0\right)},\tag{11}$$

and the predictive distribution of $X_{i,d,t}$

$$Q_{i,d}\left(dx_{i,d,t+1}\mid \mathcal{F}_t^a\right)=\frac{\int_{\theta\in\Theta}Q_{i,d}\left(dx_{i,d,t+1};\theta_i\right)L\left(\mathcal{F}_t^a;\theta\right)F\left(d\theta;\varepsilon_0\right)}{\int_{\theta\in\Theta}L\left(\mathcal{F}_t^a;\theta\right)F\left(d\theta;\varepsilon_0\right)}$$

$$=\frac{\int_{\theta\in\Theta}Q_{i,d}\left(dx_{i,d,t+1};\theta_i\right)\prod_{i=1}^{k}\prod_{l=1}^{m}\prod_{\tau=1}^{t_{i,d}}q_{i,d}\left(\tilde{X}_{i,d,\tau};\theta_{i,d}\right)F\left(d\theta;\varepsilon_0\right)}{\int_{\theta\in\Theta}\prod_{i=1}^{k}\prod_{d=1}^{m}\prod_{\tau=1}^{t_{i,d}}q_{(i,d)}\left(\tilde{X}_{i,d,\tau};\theta_{i,d}\right)F\left(d\theta;\varepsilon_0\right)}.\tag{12}$$

The (11)-(12) indicate that the posterior distribution of $\theta$ and predictive distribution of $X_{i,d,t}$ are independent of the allocation policy $A$, the information set $\mathcal{F}_T$ wholly determine the posterior distribute and the information set $\mathcal{F}_t$ wholly determine the predictive distribution. ∎

According to Theorem 1, we can formulate the SCP (7)-(9) as the following MDP with $T+1$ states:

$\mathcal{F}_t$ is the state at step $t$, $0 < t \leq T$, and $\mathcal{F}_0 = \varepsilon_0$;

$a_{t+1}(\mathcal{F}_t)$ is allocation action at step $t$, $0 < t \leq T$, and terminal selection action $S(\mathcal{F}_T)$ for $t = T$;

The state transition process of $\mathcal{F}_t$ for $0 \leq t < T$ as follows:

$$\mathcal{F}_t = \left(\varepsilon_0, X_{1,1}^{(t)},\ldots,X_{1,m}^{(t)},\ldots,X_{k,1}^{(t)},\ldots,X_{k,m}^{(t)}\right)$$
$$\Rightarrow \mathcal{F}_{t+1} \triangleq \left(\varepsilon_0, X_{1,1}^{(t)},\ldots,X_{1,m}^{(t)},\ldots,X_{i,1}^{(t)},\ldots,X_{i,d-1}^{(t)},X_{i,d}^{(t)},X_{i,d},X_{i,d+1}^{(t)},\ldots,X_{1,m}^{(t)}\right.$$
$$\left.,\ldots,X_{k,1}^{(t)},\ldots,X_{k,m}^{(t)}\right)\Big|_{(i,d)=a_{t+1}(\mathcal{F}_t)}$$
$$= \left(\varepsilon_0, X_{1,1}^{(t+1)},\ldots,X_{1,m}^{(t+1)},\ldots,X_{k,1}^{(t+1)},\ldots,X_{k,m}^{(t+)}\right)$$

where $X_{i,d,t}\sim Q_{i,d}(\cdot\mid\mathcal{F}_t)$. Moreover, this MDP only has a nonzero reward, which is the terminal reward $V_T(\mathcal{F}_T)$. Then, we can recursively capture the optimal allocation policy $A^*$ in (9) by solving the following Bellman equation:

$$V_T(\mathcal{F}_T) \triangleq \mathrm{P}\left\{\bigcap_{i=1,i\neq 1}^{k}\bigcap_{l=1}^{m}\left(\mu_{\langle 1\rangle_T,l}\geq \mu_{\langle i\rangle_T,d_{\langle i\rangle_T}}\right)\Bigg|\mathcal{F}_T\right\},\tag{13}$$

and

$$V_t(\mathcal{F}_t) \triangleq \mathrm{E}\left[V_{t+1}(\mathcal{F}_{t+1}) \mid \mathcal{F}_t\right]\Big|_{(i,d)=a^*_{t+1}(\mathcal{F}_t)}, \tag{14}$$

where $0 \leq t < T$, and

$$a^*_{t+1}(\mathcal{F}_t) = \underset{(i,d) \in \{1,2,\ldots,k\}\times\{1,2,\ldots,m\}}{\arg\max} \mathrm{E}[V_{t+1}(\mathcal{F}_{t+1}) \mid \mathcal{F}_t]. \tag{15}$$

**Remark.** According to Proposition 1.3.1 of Bertsekas (2012), the optimal allocation policy of SCP (7)-(8) is equivalent to the optimal solution of Bellman equation (10)-(12).

Since no specific form of distribution is specified for $Q(\cdot\,;\theta)$, the above framework applies to the general distribution case. To make the proposed simulation budget allocation policy suitable for practical implementation, this paper focuses on the case where $X_{i,d,t} \sim N(\mu_{i,l}, \sigma_{i,l}^2)$, $(i,d) \in \{1,2,\ldots,k\}\times\{1,2,\ldots,m\}$, $0 \leq t < T$, with unknown mean $\mu_{i,l}$ and known variance $\sigma_{i,d}^2$. By the DeGroot (2005), we have $\mu_{i,d}$ follows the conjugate Gaussian prior distribution $N\left(\mu_{i,d}^{(0)}, \left(\sigma_{i,d}^{(0)}\right)^2\right)$ and the posterior distribution of $\mu_{i,d}$ is $N\left(\mu_{i,d}^{(t)}, \left(\sigma_{i,d}^{(t)}\right)^2\right)$, where

$$\begin{aligned}
\mu_{i,d}^{(t)} &= \left(\sigma_{i,d}^{(t)}\right)^2 \left(\frac{\mu_{i,d}^{(0)}}{\left(\sigma_{i,d}^{(0)}\right)^2} + \frac{t_{i,d}\bar{X}_{i,d}^{(t)}}{\sigma_{i,d}^2}\right), \\
\left(\sigma_{i,d}^{(t)}\right)^2 &= \left(\frac{1}{\left(\sigma_{i,d}^{(0)}\right)^2} + \frac{t_{i,d}}{\sigma_{i,d}^2}\right)^{-1},
\end{aligned} \tag{16}$$

and $\bar{X}_{i,d}^{(t)} = \sum_{w=1}^{t_i} X_{i,d,w}/t_i$, $(i,d) \in \{1,2,\ldots,k\}\times\{1,2,\ldots,m\}, 0 \leq t < T$. In addition, if $\sigma_{i,d}^{(0)} \to \infty$, $\mu_{i,d}^{(t)} = \bar{X}_{i,d}^{(t)}$, such a prior is uninformative. For a Gaussian distribution with unknown mean and variance, which has a normal gamma conjugate prior. If the Gaussian assumption is not satisfied, the macroscopic replication obtained by batch means or average performances follows an approximately Gaussian distribution according to the central limit theorem (Chen & Lee 2011).

## 3. Efficient Dynamic Allocation Policy

This section follows the approximate dynamic programming paradigm to derive an efficient asymptotically optimal dynamic allocation policy with an analytical form based on a single feature value function approximation method. Specifically, to alleviate the common curse-of-dimensional problem of MDP, similar to the approach of Peng et al. (2018), this paper abandons the backward induction method and adopts forward programming by maximizing the value function approximation (VFA) one step ahead, thus finally solving the stochastic dynamic programming problem (12).

We suppose any step $t$ could be the last step. Based on the information set $\mathcal{F}_t$, we have the $\mu_{i,d}$ follows a Gaussian posterior distribution with mean $\mu_{i,d}^{(t)}$ and $\left(\sigma_{i,d}^{(t)}\right)^2$, i.e., $\mu_{i,d} \sim N\left(\mu_{i,d}^{(t)}, \left(\sigma_{i,d}^{(t)}\right)^2\right)$, where $(i,d) \in \{1,2,\ldots,k\} \times \{1,2,\ldots,m\}$. So, the vector

$$\left(\mu_{\langle 1\rangle_t,1} - \mu_{\langle 2\rangle_t,d_{\langle 2\rangle_t}}, \ldots, \mu_{\langle 1\rangle_t,1} - \mu_{\langle k\rangle_t,d_{\langle k\rangle_t}}, \ldots, \mu_{\langle 1\rangle_t,m} - \mu_{\langle 2\rangle_t,d_{\langle 2\rangle_t}}, \ldots, \mu_{\langle 1\rangle_t,m} - \mu_{\langle k\rangle_t,d_{\langle k\rangle_t}}\right)$$

follows a joint multivariate Gaussian distribution $MN(\mu^{(t)}, \Sigma)$, where

$$\mu^{(t)} = \left(\mu_{\langle 1\rangle_t,1}^{(t)} - \mu_{\langle 2\rangle_t,d_{\langle 2\rangle_t}}^{(t)}, \ldots, \mu_{\langle 1\rangle_t,1}^{(t)} - \mu_{\langle k\rangle_t,d_{\langle k\rangle_t}}^{(t)}, \ldots, \mu_{\langle 1\rangle_t,m}^{(t)} - \mu_{\langle 2\rangle_t,d_{\langle 2\rangle_t}}^{(t)}, \ldots, \mu_{\langle 1\rangle_t,m}^{(t)} - \mu_{\langle k\rangle_t,d_{\langle k\rangle_t}}^{(t)}\right)$$

and covariance matrix $\Sigma = \Gamma' \Lambda \Gamma$ is the covariance matrix, where the $'$ denotes transpose operation of the matrix, the diagonal matrix:

$$\Lambda \triangleq \mathrm{diag}\left((\sigma_{\langle 1\rangle_t,1}^{(t)})^2, (\sigma_{\langle 1\rangle_t,2}^{(t)})^2, \ldots, (\sigma_{\langle 1\rangle_t,m}^{(t)})^2, (\sigma_{\langle 2\rangle_t,d_{\langle 2\rangle_t}}^{(t)})^2, \ldots, (\sigma_{\langle k\rangle_t,d_{\langle k\rangle_t}}^{(t)})^2\right)_{(m+k-1)\times(m+k-1)},$$

and $\Gamma$ is block matrix,

$$\Gamma \triangleq \begin{pmatrix} P_1 & P_2 & \cdots & P_m \\ Q_1 & Q_2 & \cdots & Q_m \end{pmatrix}_{(m+k-1)\times(k-1)m},$$

in which $P_i \in \mathbb{R}^{m \times (k-1)(m-1)}$, $i \in \{1,\ldots,m\}$ is a matrix where $i$-th row is 1 and the remaining rows are 0, and $Q_i \in \mathbb{R}^{(k-1)\times(k-1)}$, $i \in \{1,\ldots,m\}$, is a diagonal matrix, i.e.,

$$P_i \triangleq \begin{pmatrix} 0 & 0 & \cdots & 0 \\ 0 & 0 & \cdots & 0 \\ \vdots & \vdots & \ddots & \vdots \\ 1 & 1 & \cdots & 1 \\ \vdots & \vdots & \ddots & \vdots \\ 0 & 0 & \cdots & 0 \end{pmatrix}_{m \times (k-1)}, \quad Q_i \triangleq \begin{pmatrix} -1 & 0 & \cdots & 0 \\ 0 & -1 & \cdots & 0 \\ \vdots & \vdots & \ddots & \vdots \\ 0 & 0 & \cdots & -1 \end{pmatrix}_{k-1 \times k-1}.$$

Following Cholesky decomposition, there exists an upper triangular matrix $U \triangleq (u_{i,j})_{m(k-1)\times m(k-1)}$ such that $\Sigma = U'U = (a_{i,j})_{m(k-1)\times m(k-1)}$ in which

$$\begin{cases} a_{i,j} = u_{1,i}u_{1,j} + u_{2,i}u_{2,j} + \ldots + u_{j,i}u_{j,j}, & i < j, \\ a_{i,i} = u_{1,i}^2 + u_{2,i}^2 + \ldots + u_{i,i}^2. \end{cases}$$

The value function for $S(\mathcal{F}_t) = \langle 1 \rangle_t$ can be rewritten

$$V_t(\mathcal{F}_t) = P\left\{\bigcap_{i=1,i\neq 1}^{k}\bigcap_{l=1}^{m}\left(\mu_{\langle 1\rangle_t,l} \geq \mu_{\langle i\rangle_t,d_{\langle i\rangle_t}}\right)\bigg|\mathcal{F}_t\right\}$$
$$= P\left\{\mu_{\langle 1\rangle_t,1} - \mu_{\langle 2\rangle_t,d_{\langle 2\rangle_t}} > 0, \ldots, \mu_{\langle 1\rangle_t,l} - \mu_{\langle k\rangle_t,d_{\langle k\rangle_t}} > 0\right.$$
$$\left.,\ldots, \mu_{\langle 1\rangle_t,m} - \mu_{\langle 2\rangle_t,d_{\langle 2\rangle_t}} > 0, \ldots, \mu_{\langle 1\rangle_t,m} - \mu_{\langle k\rangle_t,d_{\langle k\rangle_t}} > 0\right\} \quad (17)$$
$$= \frac{1}{(\sqrt{2\pi})^{m(k-1)}}\iint_Z \prod_{j=1}^{m(k-1)} \exp(-\frac{z_j^2}{2})d_{z_1},\ldots,d_{z_{m(k-1)}},$$

where $Z_j, j \in \{1,\ldots,m(k-1)\}$ are independent standard Gaussian random variables and the integration region $Z$ is enclosed by some hyperplanes:

$$\sum_{\gamma=1}^{q} u_{\gamma,q} z_\gamma > \mu_{\langle i\rangle_t,d_{\langle i\rangle_t}}^{(t)} - \mu_{\langle 1\rangle_t,l}^{(t)}, \quad q=1,\ldots,m(k-1). \quad (18)$$

where $i \in \{2,\ldots,k\}, l \in \{1,\ldots,m\}$.

Since the density function of the multivariate standard Gaussian distribution decreases at an exponential rate with respect to the distance from the origin, the hypersphere size centered at the origin can capture the integral value over the entire integration region $Z$, i.e., the value function (17) (Peng et al., 2018). The advantage of approximating (17) by the size of the inner tangent hypersphere is that the sequence of hypersphere sizes is equivalent to the sequence of sizes of its radius. Therefore, the problem of maximizing the value function (17) can be equivalently transformed into the above problem of maximizing the radius of the inner tangent hypersphere, and the approximation error decreases to zero at an exponential rate when the radius of the maximum inner tangent sphere becomes infinite.

Based on the above analysis, the remainder of this work will focus on solving the problem of maximizing the above inner tangent hypersphere radius. we use the square of the radius of the largest inner tangent hypersphere in the integration region $Z$ as a single feature value function approximation of the value function (14) (Powell, 2007), i.e.,

$$\tilde{V}_t(\mathcal{F}_t) = \min\left\{R_{1,1,2,d_2}^2(\mathcal{F}_t),\ldots,R_{1,1,k,d_k}^2(\mathcal{F}_t),\ldots,R_{1,m,2,d_2}^2(\mathcal{F}_t),\ldots,R_{1,m,k,d_k}^2(\mathcal{F}_t)\right\}, \quad (19)$$

where

$$R_{1,l,i,d_i}^2(\mathcal{F}_t) = \frac{\left(\mu_{\langle 1\rangle_t,l}^{(t)} - \mu_{\langle i\rangle_t,d_{\langle i\rangle_t}}^{(t)}\right)^2}{\left(\left(\sigma_{\langle 1\rangle_t,l}^{(t)}\right)^2 + \left(\sigma_{\langle i\rangle_t,d_{\langle i\rangle_t}}^{(t)}\right)^2\right)^2},$$

in which $d \in \{1,\ldots,m\}, i \in \{1,\ldots,k\}$.

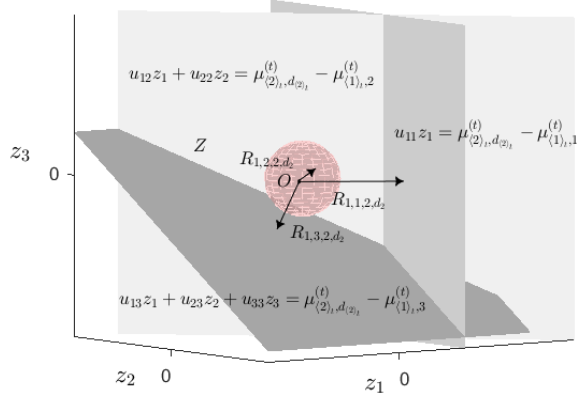

Fig. 1. Using the integrand function $\exp\left\{\frac{-(z_1^2+z_2^2+z_3^2)}{2}\right\}$ in hypersphere region to approximate its integral on the entire integral region $Z$.

In Figure 1, we show the VFA process of the value function $V_t(\mathcal{F}_t)$ with $k=2, m=3$ as an example. The (14) is the integral value of the standard Gaussian distribution over the region $Z$, and the $\tilde{V}_t(\mathcal{F}_t)$ is $R_{1,2,2d_2}$.

Suppose $t+1$-th simulation replication is the last one. A VFA of $\tilde{V}_t(\mathcal{F}_t)$ looking one step ahead at step $t$ by alternative-scenario pair $(j,s) \in \{1,2,\ldots,k\} \times \{1,2,\ldots,m\}$ as shown follows.

$$V_t(\mathcal{F}_t;(j,s)) = \mathrm{E}\left[V_{t+1}(\mathcal{F}_t, X_{j,s,t+1}) \big| \mathcal{F}_t\right]. \tag{20}$$

Following the certainty equivalent approximation in Bertsekas (2005), we make a further approximation for (20)

$$V_t(\mathcal{F}_t;(j,s)) = V_{t+1}\left(\mathcal{F}_t, \mathrm{E}[X_{j,s,t+1} | \mathcal{F}_t]\right), \tag{21}$$

where $j=1, s \in \{1,2,\ldots,m\}$, then

$$V_t(\mathcal{F}_t;(1,s)) = \min\left\{\min_{i \in \{1,\ldots,k\}, i \neq 1}\left\{\frac{\left(\mu^{(t)}_{\langle 1 \rangle_t,s} - \mu^{(t)}_{\langle i \rangle_t, d_{\langle i \rangle_t}}\right)^2}{\left(\sigma^{(t+1)}_{\langle 1 \rangle_t,s}\right)^2 + \left(\sigma^{(t)}_{\langle i \rangle_t, d_{\langle i \rangle_t}}\right)^2}\right\},\right.$$

$$\left.\min_{l \in \{1,\ldots,m\}, l \neq s, i \in \{1,\ldots,k\}, i \neq 1}\left\{\frac{\left(\mu^{(t)}_{\langle 1 \rangle_t,l} - \mu^{(t)}_{\langle i \rangle_t, d_{\langle i \rangle_t}}\right)^2}{\left(\sigma^{(t)}_{\langle 1 \rangle_t,l}\right)^2 + \left(\sigma^{(t)}_{\langle i \rangle_t, d_{\langle i \rangle_t}}\right)^2}\right\}\right\}, \tag{22}$$

and $j \in \{1,2,\ldots,k\}, j \neq 1, s \in \{d_2,\ldots,d_k\}$, then

$$V_t(\mathcal{F}_t;(j,s)) = \min\left\{\min_{l\in\{1,\ldots,m\}}\left\{\frac{\left(\mu_{\langle 1\rangle_t,l}^{(t)}-\mu_{\langle j\rangle_t,s}^{(t)}\right)^2}{\left(\sigma_{\langle 1\rangle_t,l}^{(t)}\right)^2+\left(\sigma_{\langle j\rangle_t,s}^{(t+1)}\right)^2}\right\},\right.$$
$$\left.\min_{l\in\{1,\ldots,m\},i\in\{1,\ldots,k\},i\neq 1}\left\{\frac{\left(\mu_{\langle 1\rangle_t,l}^{(t)}-\mu_{\langle i\rangle_t,d_{\langle i\rangle_t}}^{(t)}\right)^2}{\left(\sigma_{\langle 1\rangle_t,l}^{(t)}\right)^2+\left(\sigma_{\langle i\rangle_t,d_{\langle i\rangle_t}}^{(t)}\right)^2}\right\}\right\}. \quad (23)$$

An asymptotically optimal dynamic allocation policy (RAODA) that maximizes the VFA is given by

$$\hat{a}_{t+1}(\mathcal{F}_t) = \underset{(j,s)\in\{1,2,\ldots,k\}\times\{1,2,\ldots,m\}}{\operatorname{argmax}} V_t(\mathcal{F}_t;(j,s)). \quad (24)$$

**Remark**. It is worth that the problem reduces to select the best alternative with determined scenario, if $m = 1$. In this case, (19)-(20) can be simplified to the following form: for $j = 1$, then

$$V_t(\mathcal{F}_t;1) = \min_{i\in\{1,\ldots,k\},i\neq 1}\left\{\frac{\left(\mu_{\langle 1\rangle_t}^{(t)}-\mu_{\langle i\rangle_t}^{(t)}\right)^2}{\left(\sigma_{\langle 1\rangle_t}^{(t+1)}\right)^2+\left(\sigma_{\langle i\rangle_t}^{(t)}\right)^2}\right\}, \quad (25)$$

and $j \in \{1,\ldots,k\}, j \neq 1$, then

$$V_t(\mathcal{F}_t;j) = \min\left\{\frac{\left(\mu_{\langle 1\rangle_t,l}^{(t)}-\mu_{\langle j\rangle_t,s}^{(t)}\right)^2}{\left(\sigma_{\langle 1\rangle_t,l}^{(t)}\right)^2+\left(\sigma_{\langle j\rangle_t,s}^{(t+1)}\right)^2},\min_{i\in\{1,\ldots,k\},i\neq j\neq 1}\left\{\frac{\left(\mu_{\langle 1\rangle_t}^{(t)}-\mu_{\langle i\rangle_t}^{(t)}\right)^2}{\left(\sigma_{\langle 1\rangle_t}^{(t)}\right)^2+\left(\sigma_{\langle i\rangle_t}^{(t)}\right)^2}\right\}\right\} \quad (26)$$

which is the same as the simulation budget dynamic allocation policy proposed in Peng et.al (2018).

## 4. Consistency and Asymptotically Optimality

In this section, we analyze that under the allocation policy RAODA (22)-(24), the selected best and worst subsets are consistent with the truly best and worst subsets as $t \to \infty$. Moreover, the simulation budget allocation ratio obtained by the RAOD can asymptotically converge to the optimal simulation budget allocation rate proposed by Gao et al. (2017) based on OCBA method.

**Theorem 2.** As $t \to \infty$, then

$$\lim_{t\to\infty}\langle 1\rangle_t = \langle 1\rangle, \ a.s. \quad (27)$$

*Proof*: Define $\Omega = \{(1,1),\ldots,(1,m)(2,d_2),\ldots,(k,d_k)\}$. For an alternative, as long as all the other $k - 1$ alternatives have at least one scenario inferior to all the $m$ scenarios of this alternative, it is sufficient to conclude that this alternative is the best alternative. Then, we just need to show that, following the RAODA, every alternative-scenario pair $(j, s) \in \Omega$ will get an infinite number of simulation replications often almost surely, i.e., $\lim_{t \to \infty} t_{j,s} = \infty$, and the consistency will follow by the law of large numbers (LLN). Define

$$\Xi = \{(j,s) \in \Omega | \text{ alternative-scenario pair } (j,s) \text{ is allocated infinitely simulation replications often } a.s.\}.$$

We firstly prove that $\{(1,1), \ldots, (1, m)\} \cap \Xi \neq \emptyset$, otherwise, there are $(j, s) \in \{(1,1), \ldots, (1, m)\}$ and $(j, s) \in \{(2, d_2), \ldots, (k, d_k)\} \cap \Xi$ such that $\lim_{t \to \infty} t_{j,s} < \infty$, $\lim_{t \to \infty} t_{i,l} = \infty$, and

$$\lim_{t \to \infty} V_t(\mathcal{F}_t; (j,s)) - V_t(\mathcal{F}_t) > 0, . \tag{28}$$

$$\lim_{t \to \infty} V_t(\mathcal{F}_t; (i,l)) - V_t(\mathcal{F}_t) = 0, . \tag{29}$$

this contradicts with the RAODA. So, we have $\{(1,1), \ldots, (1, m)\} \cap \Xi \neq \emptyset$.

We next prove that $\{(2, d_2), \ldots, (k, d_k)\} \cap \Xi \neq \emptyset$, otherwise, there are $(j, s) \in \{(1,1), \ldots, (1, m)\}$ and $(i, l) \in \{(2, d_2), \ldots, (k, d_k)\}$ such that $\lim_{t \to \infty} t_{j,s} = \infty$, $\lim_{t \to \infty} t_{i,l} < \infty$, and

$$\lim_{t \to \infty} V_t(\mathcal{F}_t; (j,s)) - V_t(\mathcal{F}_t) = 0, . \tag{30}$$

$$\lim_{t \to \infty} V_t(\mathcal{F}_t; (i,l)) - V_t(\mathcal{F}_t) > 0, . \tag{31}$$

this contradicts with the RAODA. So, we have $\{(2, d_2), \ldots, (k, d_k)\} \cap \Xi \neq \emptyset$.

Moreover, we clam that $\{(2, d_2), \ldots, (k, d_k)\} \subset \Xi$, otherwise, there are $(j, s) \in \{(1,1), \ldots, (1, m)\}$ and $(i, l) \in \{(2, d_2), \ldots, (k, d_k)\} \cap \Xi$ such that $\lim_{t \to \infty} t_{j,s} = \infty$, $\lim_{t \to \infty} t_{i,l} < \infty$, and

$$\lim_{t \to \infty} V_t(\mathcal{F}_t; (j,s)) - V_t(\mathcal{F}_t) = 0, \tag{32}$$

$$\lim_{t \to \infty} V_t(\mathcal{F}_t; (j,s)) - V_t(\mathcal{F}_t) > 0, \tag{33}$$

this contradicts with the RAODA. So, we have $\{(2, d_2), \ldots, (k, d_k)\} \subset \Xi$.

Finally, suppose $\exists (i, l)\{(2, d_2), \ldots, (k, d_k)\}$ is only allocated finite simulation replication, then $\exists (j, s) \in \{(1,1), \ldots, (1, m)\}$ such that

$$\lim_{t \to \infty} V_t(\mathcal{F}_t; (j,s)) - V_t(\mathcal{F}_t) = 0, \tag{34}$$

$$\lim_{t\to\infty} V_t(\mathcal{F}_t;(j,s)) - V_t(\mathcal{F}_t) = 0. \tag{35}$$

So, we have $\{(1,1),\dots,(1,m)\} \subset \Xi$. Synthesize the above analysis, the RAODA can achieve consistency. ∎

**Remark.** The proof process of Theorem 2 indicates that $\lim_{t\to\infty} t_{j,s} = \infty$, $\forall (j,s) \in \Omega$. Then by the LLN, we have $\lim_{t\to\infty} \mu_{j,s}^{(t)} = \mu_{j,s}$, $\lim_{t\to\infty} \left(\sigma_{j,s}^{(t)}\right)^2 = (\sigma_{j,s})^2/t_{j,s}$, $\forall (j,s) \in \Omega$. In order to facilitate the following analysis, we use the $\mu_{j,s}^{(t)}$ and $\left(\sigma_{j,s}^{(t)}\right)^2$ instead of $\mu_{\forall(j,s)}$ and $(\sigma_{j,s})^2/t_{j,s}$ in the RAODAP, then $\hat{V}_t(\mathcal{F}_t;(j,s))$, $\forall (j,s) \in \Omega$, can be rewritten as follows

if $j = 1, s \in \{1,2,\dots,m-1\}$, then

$$V_t(\mathcal{F}_t;(1,s)) = \min\left\{\min_{i\in\{1,\dots,k\},i\neq 1}\left\{\frac{\left(\mu_{\langle 1\rangle,s}-\mu_{\langle i\rangle,d_{\langle i\rangle}}\right)^2}{\sigma^2_{\langle 1\rangle,s}/(t_{\langle 1\rangle,s}+1)+\sigma^2_{\langle i\rangle,d_{\langle i\rangle}}/t_{\langle i\rangle,d_{\langle i\rangle}}}\right\},\right.$$
$$\left.\min_{l\in\{1,\dots,m\},l\neq l'}\min_{i\in\{1,\dots,k\},i\neq 1}\left\{\frac{\left(\mu_{\langle i\rangle,s}-\mu_{\langle i\rangle,d_{\langle i\rangle}}\right)^2}{\sigma^2_{\langle i\rangle,s}/t_{\langle i\rangle,s}+\sigma^2_{\langle i\rangle,d_{\langle i\rangle}}/t_{\langle i\rangle,d_{\langle i\rangle}}}\right\}\right\} \tag{36}$$

if $j \neq 1, s \in \{d_2,\dots,d_k\}$, then

$$V_t(\mathcal{F}_t;(j,s)) = \min\left\{\min_{l\in\{1,\dots,m\}}\left\{\frac{\left(\mu_{\langle 1\rangle,l}-\mu_{\langle j\rangle,s}\right)^2}{\sigma^2_{\langle 1\rangle,l}/t_{\langle 1\rangle,l}+\sigma^2_{\langle j\rangle,s}/(t_{\langle i\rangle,s}+1)}\right\},\right.$$
$$\left.\min_{l\in\{1,\dots,m\},i\in\{1,\dots,k\},i\neq j\neq 1}\left\{\frac{\left(\mu_{\langle 1\rangle,l}-\mu_{\langle i\rangle,d_{\langle i\rangle}}\right)^2}{\sigma^2_{\langle 1\rangle,s}/t_{\langle i\rangle,l}+\sigma^2_{\langle i\rangle,d_{\langle i\rangle}}/t_{\langle i\rangle,d_{\langle i\rangle}}}\right\}\right\} \tag{37}$$

The following Theorem 3 establishes that each clustering point of the simulation budget allocation ratio $\left(\alpha_{1,1}^{(t)},\dots,\alpha_{1,m}^{(t)},\dots,\alpha_{k,1}^{(t)},\dots,\alpha_{k,m}^{(t)}\right)$ obtained by the RAODA satisfies the optimality condition for the simulation budget allocation proposed by Gao et al. (2017) based on optimal computing budget allocation method.

**Theorem 3.** Define the optimal simulation budget allocation ratio $(\alpha_{1,1}^*,\dots,\alpha_{1,m}^*,\dots,\alpha_{k,1}^*,\dots,\alpha_{k,m}^*)$, which satisfies the optimality condition as follows

$$\min_{i\in\{1,\dots,k\},i\neq 1} G_{1,l,i,d_i}\left(\alpha_{\langle i\rangle,l}^*,\alpha_{\langle i\rangle,d_{\langle i\rangle}}^*\right) = \min_{i\in\{1,\dots,k\},i\neq 1} G_{1,l',i,d_i}\left(\alpha_{\langle 1\rangle,l'}^*,\alpha_{\langle i\rangle,d_{\langle i\rangle}}^*\right),\ l,l'\in\{1,..,m\},l\neq l', \tag{38}$$

$$\min_{l\in\{1,\dots,m\}} G_{1,l,i,d_i}\left(\alpha_{\langle i\rangle,l}^*,\alpha_{\langle i\rangle,d_{\langle i\rangle}}^*\right) = \min_{l\in\{1,\dots,m\}} G_{1,l,i',d_i}\left(\alpha_{\langle 1\rangle,l}^*,\alpha_{\langle i'\rangle,d_{\langle i\rangle}}^*\right),\ i,i'\in\{1,..,k\},i\neq i'\neq 1, \tag{39}$$

$$\sum_{l=1}^{m} \frac{\left(\alpha_{\langle 1\rangle,l}^{*}\right)^{2}}{\sigma_{\langle 1\rangle,l}^{2}} = \sum_{i=2}^{k} \frac{\left(\alpha_{\langle i\rangle,d_{\langle i\rangle}}^{*}\right)^{2}}{\sigma_{\langle i\rangle,d_{\langle i\rangle}}^{2}}, \tag{40}$$

$$\sum_{j=1}^{k}\sum_{s=1}^{m} \alpha_{j,s}^{*} = 1, \alpha_{j,s}^{*} \geq 0, \forall (j,s) \in \{1,2,...,k\} \times \{1,2,...,m\}, \tag{41}$$

where

$$G_{1,l,i,d_i}\left(\alpha_{\langle 1\rangle,l}, \alpha_{\langle i\rangle,d_i}\right) = \frac{\left(\mu_{\langle 1\rangle,l} - \mu_{\langle i\rangle,d_i}\right)^2}{2((\sigma_{\langle 1\rangle,l})^2 / \alpha_{\langle 1\rangle,l} + (\sigma_{\langle i\rangle,d_i})^2 / \alpha_{\langle i\rangle,d_i})}. \tag{42}$$

Then, the clustering points of $\left(\alpha_{1,1}^{(t)}, ..., \alpha_{1,m}^{(t)}, ..., \alpha_{k,1}^{(t)}, ..., \alpha_{k,m}^{(t)}\right)$ satisfy the optimality condition (38)-(40).

*Proof:* By Bolzano-Weierstrass theorem (Rudin 1976), there exist a subsequence of $\left(\alpha_{1,1}^{(t)}, ..., \alpha_{1,m}^{(t)} \alpha_{2,d_2}^{(t)}, ..., \alpha_{k,d_k}^{(t)}\right)$ converges to $(\hat{\alpha}_{1,1}, ..., \hat{\alpha}_{1,m}, \hat{\alpha}_{2,d_2}, ..., \hat{\alpha}_{k,d_k})$ such that $\hat{\alpha}_{j,s} > 0, (j,s) = \Omega$. Without loss of generality, we can suppose that $\left(\alpha_{1,1}^{(t)}, ..., \alpha_{1,m}^{(t)} \alpha_{2,d_2}^{(t)}, ..., \alpha_{k,d_k}^{(t)}\right)$ converges to $(\hat{\alpha}_{1,1}, ..., \hat{\alpha}_{1,m}, \hat{\alpha}_{2,d_2}, ..., \hat{\alpha}_{k,d_k})$. By noticing that

$$\lim_{t \to +\infty}\left[ \frac{\left(\mu_{\langle 1\rangle,l} - \mu_{\langle j\rangle,d_j}\right)^2}{\sigma_{\langle 1\rangle,l}^2 / t_{\langle i\rangle,l} + \sigma_{\langle j\rangle,s}^2 / \left(t_{\langle j\rangle,d_j} + 1\right)} - \frac{\left(\mu_{\langle 1\rangle,l} - \mu_{\langle i\rangle,d_i}\right)^2}{\sigma_{\langle 1\rangle,l}^2 / t_{\langle 1\rangle,l} + \sigma_{\langle i\rangle,d_i}^2 / t_{\langle i\rangle,d_i}} \right]$$

$$= \lim_{t \to +\infty} t \left[ \frac{\left(\mu_{\langle 1\rangle,l} - \mu_{\langle i\rangle,d_i}\right)^2}{\sigma_{\langle 1\rangle,l}^2 / \alpha_{\langle 1\rangle,l} + \sigma_{\langle i\rangle,d_i}^2 / \left(\alpha_{\langle i\rangle,d_i} + 1/t\right)} - \frac{\left(\mu_{\langle 1\rangle,l} - \mu_{\langle i\rangle,d_i}\right)^2}{\sigma_{\langle 1\rangle,l}^2 / \alpha_{\langle 1\rangle,l} + \sigma_{\langle i\rangle,d_i}^2 / \alpha_{\langle i\rangle,d_i}} \right]$$

$$= \lim_{t \to +\infty} \frac{\partial M_{1,l,i,d_i}\left(\alpha_{\langle 1\rangle,l}^{(t)}, x\right)}{\partial x} \Bigg|_{x = \alpha_{\langle i\rangle,d_i}^{(t)}} = \lim_{t \to +\infty} \left(\frac{\sigma_{\langle i\rangle,d_i}}{\alpha_{\langle i\rangle,d_i}^{(t)}}\right)^2 \frac{\left(\mu_{\langle 1\rangle,l} - \mu_{\langle i\rangle,d_i}\right)^2}{\sigma_{\langle 1\rangle,l}^2 / \alpha_{\langle 1\rangle,l} + \sigma_{\langle i\rangle,d_i}^2 / \alpha_{\langle i\rangle,d_i}}, \tag{43}$$

and

$$\lim_{t \to +\infty}\left[ \frac{\left(\mu_{\langle 1\rangle,l} - \mu_{\langle j\rangle,d_j}\right)^2}{\sigma_{\langle 1\rangle,l}^2 / \left(t_{\langle i\rangle,l} + 1\right) + \sigma_{\langle j\rangle,s}^2 / t_{\langle j\rangle,d_j}} - \frac{\left(\mu_{\langle 1\rangle,l} - \mu_{\langle i\rangle,d_i}\right)^2}{\sigma_{\langle 1\rangle,l}^2 / t_{\langle 1\rangle,l} + \sigma_{\langle i\rangle,d_i}^2 / t_{\langle i\rangle,d_i}} \right]$$

$$= \lim_{t \to +\infty} t \left[ \frac{\left(\mu_{\langle 1\rangle,l} - \mu_{\langle i\rangle,d_i}\right)^2}{\sigma_{\langle 1\rangle,l}^2 / \left(\alpha_{\langle 1\rangle,l} + 1\right) + \sigma_{\langle i\rangle,d_i}^2 / \alpha_{\langle i\rangle,d_i}} - \frac{\left(\mu_{\langle 1\rangle,l} - \mu_{\langle i\rangle,d_i}\right)^2}{\sigma_{\langle 1\rangle,l}^2 / \alpha_{\langle 1\rangle,l} + \sigma_{\langle i\rangle,d_i}^2 / \alpha_{\langle i\rangle,d_i}} \right]$$

$$= \lim_{t \to +\infty} \frac{\partial M_{1,l,i,d_i}\left(x, \alpha_{\langle i\rangle,d_i}^{(t)}\right)}{\partial x} \Bigg|_{x = \alpha_{\langle 1\rangle,l}^{(t)}} = \lim_{t \to +\infty} \left(\frac{\sigma_{\langle 1\rangle,l}}{\alpha_{\langle 1\rangle,l}^{(t)}}\right)^2 \frac{\left(\mu_{\langle 1\rangle,l} - \mu_{\langle i\rangle,d_i}\right)^2}{\sigma_{\langle 1\rangle,l}^2 / \alpha_{\langle 1\rangle,l} + \sigma_{\langle i\rangle,d_i}^2 / \alpha_{\langle i\rangle,d_i}}, \tag{44}$$

in which

$$M_{1,l,i,d_i}\left(\alpha_{\langle 1\rangle,l}^{(t)},\alpha_{\langle i\rangle,d_i}^{(t)}\right) \triangleq 2G_{1,l,i,d_i}\left(\alpha_{\langle 1\rangle,l}^{(t)},\alpha_{\langle i\rangle,d_i}^{(t)}\right) = \frac{\left(\mu_{\langle 1\rangle,l}-\mu_{\langle i\rangle,d_i}\right)^2}{\sigma_{\langle 1\rangle,l}^2/t_{\langle 1\rangle,l}+\sigma_{\langle i\rangle,d_i}^2/t_{\langle i\rangle,d_i}}, \tag{45}$$

$l \in \{1,2,\ldots,m\}$ and $i \in \{2,\ldots,k\}$.

We firstly prove $\hat{\alpha}_{j,s} > 0, \forall (j,s) \in \Omega$ in four steps, as follows:

Step1: We clam that $\exists (1,l) \in \{(1,1),\ldots,(1,m)\}$, such that $\hat{\alpha}_{1,l} > 0$, otherwise, there exists $(i,d_i) \in \{(2,d_2),\ldots,(k,d_k)\}$, such that $\hat{\alpha}_{j,s} > 0$, and

$$\begin{aligned}
\lim_{t\to+\infty}\left(\frac{\sigma_{\langle i\rangle,d_i}}{\alpha_{\langle i\rangle,d_i}^{(t)}}\right)^2 \frac{\left(\mu_{\langle 1\rangle,l}-\mu_{\langle i\rangle,d_i}\right)^2}{\sigma_{\langle 1\rangle,l}^2/\alpha_{\langle 1\rangle,l}+\sigma_{\langle i\rangle,d_i}^2/\alpha_{\langle i\rangle,d_i}} &= 0, \\
\lim_{t\to+\infty}\left(\frac{\sigma_{\langle 1\rangle,l}}{\alpha_{\langle 1\rangle,l}^{(t)}}\right)^2 \frac{\left(\mu_{\langle 1\rangle,l}-\mu_{\langle i\rangle,d_i}\right)^2}{\sigma_{\langle 1\rangle,l}^2/\alpha_{\langle 1\rangle,l}+\sigma_{\langle i\rangle,d_i}^2/\alpha_{\langle i\rangle,d_i}} &> 0,
\end{aligned} \tag{46}$$

for $(1,l) \in \{(1,1),\ldots,(1,m)\}$. (46) contradict with the RAODA, (36)-(37). So, we get $\exists (1,l) \in \{(1,1),\ldots,(1,m)\}$, such that $\hat{\alpha}_{1,l} > 0$.

Step2: We clam that $\exists (i,d_i) \in \{(2,d_2),\ldots,(k,d_k)\}$, such that $\hat{\alpha}_{j,s} > 0$, otherwise, there exists $(1,l) \in \{(1,1),\ldots,(1,m)\}$, such that $\hat{\alpha}_{1,l} > 0$, and

$$\begin{aligned}
\lim_{t\to+\infty}\left(\frac{\sigma_{\langle i\rangle,d_i}}{\alpha_{\langle i\rangle,d_i}^{(t)}}\right)^2 \frac{\left(\mu_{\langle 1\rangle,l}-\mu_{\langle i\rangle,d_i}\right)^2}{\sigma_{\langle 1\rangle,l}^2/\alpha_{\langle 1\rangle,l}+\sigma_{\langle i\rangle,d_i}^2/\alpha_{\langle i\rangle,d_i}} &> 0, \\
\lim_{t\to+\infty}\left(\frac{\sigma_{\langle 1\rangle,l}}{\alpha_{\langle 1\rangle,l}^{(t)}}\right)^2 \frac{\left(\mu_{\langle 1\rangle,l}-\mu_{\langle i\rangle,d_i}\right)^2}{\sigma_{\langle 1\rangle,l}^2/\alpha_{\langle 1\rangle,l}+\sigma_{\langle i\rangle,d_i}^2/\alpha_{\langle i\rangle,d_i}} &= 0,
\end{aligned} \tag{47}$$

for $(i,d_i) \in \{(2,d_2),\ldots,(k,d_k)\}$. (47) contradict with the RAODA, (36)-(37). So, we get $\exists (i,d_i) \in \{(2,d_2),\ldots,(k,d_k)\}$, such that $\hat{\alpha}_{j,s} > 0$.

Step3: We clam that $\forall (i,d_i) \in \{(2,d_2),\ldots,(k,d_k)\}$, such that $\hat{\alpha}_{1,l} > 0$, otherwise, there exists $\exists (i,d_i) \in \{(2,d_2),\ldots,(k,d_k)\}$ and $(1,l) \in \{(1,1),\ldots,(1,m)\}$ such that $\hat{\alpha}_{i,d_i} = 0$, $\hat{\alpha}_{1,l} > 0$, and

$$\begin{aligned}
\lim_{t\to+\infty}\left(\frac{\sigma_{\langle i\rangle,d_i}}{\alpha_{\langle i\rangle,d_i}^{(t)}}\right)^2 \frac{\left(\mu_{\langle 1\rangle,l}-\mu_{\langle i\rangle,d_i}\right)^2}{\sigma_{\langle 1\rangle,l}^2/\alpha_{\langle 1\rangle,l}+\sigma_{\langle i\rangle,d_i}^2/\alpha_{\langle i\rangle,d_i}} &> 0, \\
\lim_{t\to+\infty}\left(\frac{\sigma_{\langle 1\rangle,l}}{\alpha_{\langle 1\rangle,l}^{(t)}}\right)^2 \frac{\left(\mu_{\langle 1\rangle,l}-\mu_{\langle i\rangle,d_i}\right)^2}{\sigma_{\langle 1\rangle,l}^2/\alpha_{\langle 1\rangle,l}+\sigma_{\langle i\rangle,d_i}^2/\alpha_{\langle i\rangle,d_i}} &= 0,
\end{aligned} \tag{48}$$

(48) contradict with the RAODAP, (36)-(37). Thus, we obtain $\forall (i,d_i) \in \{(2,d_2)\ldots,(k,d_k)\}$, such that $\hat{\alpha}_{1,l} > 0$.

Step4: We clam that $\forall (1,l) \in \{(1,1),\ldots,(1,m)\}$, such that $\hat{\alpha}_{1,l} > 0$, otherwise, there exists $\exists (1,l) \in \{(1,1),\ldots,(1,m)\}$ such that $\hat{\alpha}_{1,l} = 0$ and

$$\lim_{t\to+\infty}\left(\frac{\sigma_{\langle i\rangle,d_i}}{\alpha^{(t)}_{\langle i\rangle,d_i}}\right)^2 \frac{\left(\mu_{\langle 1\rangle,l}-\mu_{\langle i\rangle,d_i}\right)^2}{\sigma^2_{\langle 1\rangle,l}/\alpha_{\langle 1\rangle,l}+\sigma^2_{\langle i\rangle,d_i}/\alpha_{\langle i\rangle,d_i}} > 0,$$

$$\lim_{t\to+\infty}\left(\frac{\sigma_{\langle 1\rangle,l}}{\alpha^{(t)}_{\langle 1\rangle,l}}\right)^2 \frac{\left(\mu_{\langle 1\rangle,l}-\mu_{\langle i\rangle,d_i}\right)^2}{\sigma^2_{\langle 1\rangle,l}/\alpha_{\langle 1\rangle,l}+\sigma^2_{\langle i\rangle,d_i}/\alpha_{\langle i\rangle,d_i}} = 0,$$
(49)

for $\forall(i,d_i) \in \{(2,d_2),\ldots,(k,d_k)\}$. (49) contradict with the RAODA, (36)-(37). Thus, we obtain $\hat{\alpha}_{j,s} > 0, \forall(j,s) \in \Omega$.

We next prove that $(\hat{\alpha}_{1,1},\ldots,\hat{\alpha}_{1,m},\hat{\alpha}_{2,d_2},\ldots,\hat{\alpha}_{k,d_k})$ satisfy (38)-(39), otherwise, there exists two cases as follows:

i: $\exists l, l' \in \{1,2,\ldots,m\}$ such that

$$\min_{i\in\{1,\ldots,k\},i\neq 1} M_{1,l,i,d_i}\left(\alpha^*_{\langle 1\rangle,l},\alpha^*_{\langle i\rangle,d_{\langle i\rangle}}\right) > \min_{i\in\{1,\ldots,k\},i\neq 1} M_{1,l',i,d_i}\left(\alpha^*_{\langle 1\rangle,l'},\alpha^*_{\langle i\rangle,d_{\langle i\rangle}}\right). \tag{50}$$

ii: $\exists i, i' \in \{2,\ldots,k\}, i \neq i'$ such that

$$\min_{l\in\{1,\ldots,m\}} M_{1,l,i,d_i}\left(\alpha^*_{\langle 1\rangle,l},\alpha^*_{\langle i\rangle,d_{\langle i\rangle}}\right) = \min_{l\in\{1,\ldots,m\}} M_{1,l,i',d_{i'}}\left(\alpha^*_{\langle 1\rangle,l},\alpha^*_{\langle i'\rangle,d_{\langle i'\rangle}}\right). \tag{51}$$

for any inequality in (50)-(51) to hold, this contradicts with $\left(\alpha^{(t)}_{1,1},\ldots,\alpha^{(t)}_{1,m}\alpha^{(t)}_{2,d_2},\ldots,\alpha^{(t)}_{k,d_k}\right)$ converges to $(\hat{\alpha}_{1,1},\ldots,\hat{\alpha}_{1,m},\hat{\alpha}_{2,d_2},\ldots,\hat{\alpha}_{k,d_k})$. Let us take case i as an example to prove this fact. If the inequality (50) holds, by the continuity of $M_{1,l,i,d_i}$ on $(0,1) \times (0,1)$, then there exists $N_0$ such that for $\forall t > N_0$, we have

$$\min_{i\in\{1,\ldots,k\},i\neq 1} M_{1,l,i,d_i}\left(\alpha^{(t)}_{\langle 1\rangle,l},\alpha^{(t)}_{\langle i\rangle,d_{\langle i\rangle}}\right) > \min_{i\in\{1,\ldots,k\},i\neq 1} M_{1,l',i,d_i}\left(\alpha^{(t)}_{\langle 1\rangle,l'},\alpha^{(t)}_{\langle i\rangle,d_{\langle i\rangle}}\right). \tag{52}$$

from the definition of $\alpha^{(t)}_{1,l}, \alpha^{(t)}_{1,l'}, \alpha^{(t)}_{i,d_i}, M_{1,l,i,d_i}, M_{1,l',i,d_i}$, we get that

$$\min_{i\in\{1,\ldots,k\},i\neq 1}\frac{\left(\mu_{\langle 1\rangle,l}-\mu_{\langle i\rangle,d_{\langle i\rangle}}\right)^2}{\sigma^2_{\langle 1\rangle,l}/t_{\langle 1\rangle,l}+\sigma^2_{\langle i\rangle,d_{\langle i\rangle}}/t_{\langle i\rangle,d_{\langle i\rangle}}} > \min_{i\in\{1,\ldots,k\},i\neq 1}\left\{\frac{\left(\mu_{\langle i\rangle,l'}-\mu_{\langle i\rangle,d_{\langle i\rangle}}\right)^2}{\sigma^2_{\langle i\rangle,l'}/t_{\langle i\rangle,l'}+\sigma^2_{\langle i\rangle,d_{\langle i\rangle}}/t_{\langle i\rangle,d_{\langle i\rangle}}}\right\}, \tag{53}$$

by the RAODA, (36)-(37), the alternative $(\langle 1\rangle, l)$ will be stop receiving simulation replications until the inequality sign above reverse, which contradicts with $\left(\alpha^{(t)}_{1,1},\ldots,\alpha^{(t)}_{1,m}\alpha^{(t)}_{2,d_2},\ldots,\alpha^{(t)}_{k,d_k}\right)$ converges to $(\hat{\alpha}_{1,1},\ldots,\hat{\alpha}_{1,m},\hat{\alpha}_{2,d_2},\ldots,\hat{\alpha}_{k,d_k})$. Summarizing the above $(\hat{\alpha}_{1,1},\ldots,\hat{\alpha}_{1,m},\hat{\alpha}_{2,d_2},\ldots,\hat{\alpha}_{k,d_k})$ must satisfy (38)-(39).

Finally, we show that $(\hat{\alpha}_{1,1},\ldots,\hat{\alpha}_{1,m},\hat{\alpha}_{2,d_2},\ldots,\hat{\alpha}_{k,d_k})$ must hold (38) as follows:

Case 1: If $m \leq k-1$, the (38)-(39), $\sum_{(j,s)\in\Omega} \hat{\alpha}_{j,s}=1$, and the implicit function theorem (Rudin 1964) determine implicit functions $\hat{\alpha}_{\langle i\rangle,d_i}(x_1,\ldots,x_m)|_{x_l=\hat{\alpha}_{\langle 1\rangle,l}}$, $l = 1,\ldots,m$, $i = 2,\ldots,k$, since

$$\det(\tilde{\Psi}) = \prod_{i=1}^{m} \psi_{l^{(i)},i} \left\{ \sum_{i=1}^{m} \psi_{l^{(i)},i}^{-1} \right\} > 0 \tag{54}$$

In which, $l^i \in \underset{l=\{1,\ldots,m\}}{\operatorname{argmin}} M_{1,l,i,d_i}(\alpha_{\langle 1\rangle,1}^{(t)}, \alpha_{\langle i\rangle,d_{\langle i\rangle}}^{(t)})$, $i \in \{2,\ldots,k\}$,

$$\psi_{l,i} \triangleq \left.\frac{\partial M_{\langle 1\rangle,l,i,d_i}(\hat{\alpha}_{\langle 1\rangle,l},y)}{\partial y}\right|_{y=\hat{\alpha}_{\langle i\rangle,d_{\langle i\rangle}}}, l \in \{1,\ldots,m\}, i \in \{2,\ldots,k\}, \tag{55}$$

And $\tilde{\Psi} \in \mathbb{R}^{(k-1)\times(k-1)}$,

$$\tilde{\Psi} \triangleq \begin{pmatrix} \psi_{l^{(2)},2} & -\psi_{l^{(3)},3} & \cdots & 0 & 0 \\ 0 & \psi_{l^{(2)},2} & \cdots & 0 & 0 \\ \vdots & \vdots & \cdots & \vdots & \vdots \\ 0 & 0 & \ddots & \psi_{l^{(k-1)},k-1} & -\psi_{l^{(k)},k} \\ 1 & 1 & \cdots & 1 & 1 \end{pmatrix}. \tag{56}$$

Taking the derivations of (45) and $\sum_{(j,s)\in\Omega} \hat{\alpha}_{j,s} = 1$ with respect to $\hat{\alpha}_{\langle 1\rangle,l}$, $l \in \{1,\ldots,m\}$, we can get that $\tilde{\Psi} \cdot \tilde{\Phi}_l = \Upsilon_l$, $l \in \{1,\ldots,m\}$, where $\tilde{\Phi}_l \in \mathbb{R}^{(k-1)\times 1}$ and $\tilde{\Upsilon}_l \in \mathbb{R}^{(k-1)\times 1}$,

$$\tilde{\Phi}_l \triangleq \left(\frac{\partial \tilde{\alpha}_{\langle 2\rangle,d_{\langle 2\rangle}}(y)}{\partial y},\ldots,\frac{\partial \tilde{\alpha}_{\langle k\rangle,d_{\langle k\rangle}}(y)}{\partial y}\right)'\bigg|_{y=\tilde{\alpha}_{\langle 1\rangle,l}}, l \in \{1,\ldots,m\}, \tag{57}$$

$$\gamma_{l,i} \triangleq \left.\frac{\partial M_{\langle 1\rangle,l,i,d_i}(y,\hat{\alpha}_{\langle i\rangle,d_{\langle i\rangle}})}{\partial y}\right|_{y=\hat{\alpha}_{\langle 1\rangle,l}}, l \in \{1,\ldots,m\}, i \in \{2,\ldots,k\}, \tag{58}$$

$$\gamma_i^{(l)} \triangleq \begin{cases} \gamma_{l,i+1}-\gamma_{l,i}, & l^{(i)}=l^{(i+1)}=l \\ -\gamma_{l,i}, & l^{(i)}=l, l^{(i+1)}\neq l \\ \gamma_{l,i+1}, & l^{(i)}\neq l, l^{(i+1)}=l \\ 0, & l^{(i)}\neq l, l^{(i+1)}\neq l \end{cases}, l=1,\ldots,m, i=2,\ldots,k-1, \tag{59}$$

$$\Upsilon_l \triangleq \left.\left(\gamma_2^{(l)},\ldots,\gamma_k^{(l)},-1\right)\right|_{y=\tilde{\alpha}_{\langle 1\rangle,l}}, l \in \{1,\ldots,m\}, \tag{60}$$

Where $'$ denotes the transpose operation of the matrix. In addition, by $\hat{\alpha}_{\langle 1\rangle,l}$ and $\hat{\alpha}_{\langle 1\rangle,l'}$, $l \neq l'$, $l, l' \in \{1,\ldots,m\}$ are independent, we have

$$\left.\frac{\partial M_{1,l,i,d_i}(y,\hat{\alpha}_{\langle i\rangle,d_i})}{\partial y}\right|_{y=\hat{\alpha}_{\langle 1\rangle,l'}} = 0, i \in \{2,\ldots,k\}, l' \neq l. \tag{61}$$

And, following the RAODA, (36)-(37), we can get that for $l \in \{1,\ldots,m\}$

$$\left.\frac{\partial M_{1,l^{(i)},i,d_i}(y,\hat{\alpha}_{\langle i\rangle,d_i})}{\partial y}\right|_{y=\hat{\alpha}_{\langle 1\rangle,l}} + \left.\frac{\partial M_{1,l^{(i)},i,d_j}(\hat{\alpha}_{\langle 1\rangle,l^{(i)}},y,)}{\partial y}\right|_{y=\langle i\rangle,d_i} \left.\frac{\partial \hat{\alpha}_{\langle i\rangle,d_i}(y)}{\partial y}\right|_{y=\hat{\alpha}_{\langle 1\rangle,l}} = 0, i\in\{2,...,k\}. \quad (62)$$

in which if $l^{(i)} \neq \langle l \rangle$, $\left.\frac{\partial M_{1,l^{(i)},i,d_i}(y,\hat{\alpha}_{\langle i\rangle,d_{\langle i\rangle}})}{\partial y}\right|_{y=\hat{\alpha}_{\langle 1\rangle,l}} = 0$. Otherwise, $\exists j \in \{2,...,k\}$ such that the equality above does not hold, i.e,

$$\left.\frac{\partial M_{1,l^{(j)},j,d_j}(y,\hat{\alpha}_{\langle j\rangle,d_j})}{\partial y}\right|_{y=\hat{\alpha}_{\langle 1\rangle,l}} + \left.\frac{\partial M_{1,l^{(j)},j,d_j}(\hat{\alpha}_{\langle 1\rangle,l^{(j)}},y,)}{\partial y}\right|_{y=\langle j\rangle,d_j} \left.\frac{\partial \hat{\alpha}_{\langle j\rangle,d_j}(y)}{\partial y}\right|_{y=\hat{\alpha}_{\langle 1\rangle,l}} > 0. \quad (63)$$

Following the RAODA, (36)-(37), alternative-scenario pair $(\langle 1\rangle, l)$ will receiving simulation replications and alternative-scenario pair $(\langle i\rangle, d_{\langle i\rangle})$ will stop being allocated simulation budget before the inequality sign above is no longer hold, which contradicts with $\left(\alpha_{1,1}^{(t)},...,\alpha_{1,m}^{(t)}\alpha_{2,d_2}^{(t)},...,\alpha_{k,d_k}^{(t)}\right)$ converges to $(\hat{\alpha}_{1,1},...,\hat{\alpha}_{1,m},\hat{\alpha}_{2,d_2},...,\hat{\alpha}_{k,d_k})$. Thus, we have $\widetilde{K} \cdot \widetilde{\Phi}_l = G_l, l \in \{1,...,m\}$, where $\widetilde{K} \in \mathbb{R}^{(k-1)\times(k-1)}$, $G_l \in \mathbb{R}^{(k-1)\times 1}$, and the diagonal matrix

$$\widetilde{K} \triangleq \mathrm{diag}\left(\psi_{l^{(2)},2},\psi_{l^{(3)},3},...,\psi_{l^{(k)},k}\right), \quad (64)$$

and

$$g_i^{(l)} \triangleq \begin{cases} -\gamma_{l,i+1} & l^{(i+1)} = \langle l\rangle \\ 0 & l^{(i+1)} \neq \langle l\rangle \end{cases}, l=1,2,...,m, i=1,2,...,k-1, \quad (65)$$

$$G_l \triangleq \left(g_1^{(l)}, g_2^{(l)},...,g_{k-1}^{(l)}\right)', \quad (66)$$

in which $'$ denotes the transpose operation of the matrix. According to the Lemma 3 in Zhang et al., (2022).

We have $\psi_{l,i} > 0, l = 1,...,m, i = 2,...,k$, which indicates $\det(\widetilde{K}) > 0$ and $\widetilde{K}$ is invertible,

$$\widetilde{K}^{-1} \triangleq \mathrm{diag}\left(\psi_{l^{(2)},2}^{-1},\psi_{l^{(3)},3}^{-1},...,\psi_{l^{(k)},k}^{-1}\right), \quad (67)$$

So, we have $\Upsilon_l = \widetilde{\Psi} \cdot \widetilde{K}^{-1} G_l, l = 1,...,m$, which leads to

$$\sum_{\{i|(\langle 1\rangle,l^{(i)})=(\langle 1\rangle,l)\}} \frac{\partial M_{1,l^{(i)},i,d_i}(y,\hat{\alpha}_{\langle i\rangle,d_i})/\partial y|_{y=\hat{\alpha}_{\langle 1\rangle,l^{(i)}}}}{\partial M_{1,l^{(i)},i,d_i}(\hat{\alpha}_{\langle 1\rangle,l^{(i)}},y)/\partial y|_{y=\hat{\alpha}_{\langle i\rangle,d_i}}} = 1, l=1,...,m. \quad (68)$$

Thus, we can get that

$$\frac{\left(\alpha_{\langle 1\rangle,l}^*\right)^2}{\sigma_{\langle 1\rangle,l}^2} = \sum_{\{i|(\langle 1\rangle,l^{(i)})=(\langle 1\rangle,l)\}} \frac{\left(\alpha_{\langle i\rangle,d_{\langle i\rangle}}^*\right)^2}{\sigma_{\langle i\rangle,d_{\langle i\rangle}}^2}, l=1,...,m. \quad (69)$$

Since $\cup_{l=1,...,m}\{i|(\langle 1\rangle,l^{(i)}) = (\langle 1\rangle,l)\} = \{2,...,k\}$ and $\cap_{l=1,...,m}\{i|(\langle 1\rangle,l^{(i)}) = (\langle 1\rangle,l)\} = \emptyset$,

we have (40) must hold.

Case 2: If $m > k - 1$, then the (38)-(39), $\sum_{(j,s)\in\Omega} \hat{\alpha}_{j,s}=1$, and the implicit function theorem (Rudin 1964) determine implicit functions $\hat{\alpha}_{\langle 1 \rangle,l}(x_{m+1},\ldots,x_{k+m-1})|_{x_{m+1}=\hat{\alpha}_{\langle i \rangle,d_i}}$, $l = 1,\ldots,m, i = 2,\ldots,k$, since

$$\det(\hat{\Psi}) = \prod_{i=1}^{m} \psi_{l,i^{(l)}} \left\{ \sum_{i=1}^{m} \psi_{l,i^{(l)}}^{-1} \right\} > 0 \tag{70}$$

In which, $i^{(l)} \in \underset{i=\{2,\ldots,k\}}{\arg\min} M_{1,l,i,d_i}(\alpha_{\langle 1 \rangle,1}^{(t)}, \alpha_{\langle i \rangle,d_{\langle i \rangle}}^{(t)}), l \in \{1,\ldots,m\}$,

$$\psi_{l,i} \triangleq \left. \frac{\partial M_{1,l,i,d_i}(\hat{\alpha}_{\langle 1 \rangle,l}, y)}{\partial y} \right|_{y=\hat{\alpha}_{\langle i \rangle, d_{\langle i \rangle}}}, l \in \{1,\ldots,m\}, i \in \{2,\ldots,k\}, \tag{71}$$

And $\hat{\Psi} \in \mathbb{R}^{(k-1)\times(k-1)}$,

$$\tilde{\Psi} \triangleq \begin{pmatrix} \psi_{1,i^{(1)}} & -\psi_{2,i^{(2)}} & \cdots & 0 & 0 \\ 0 & \psi_{2,i^{(2)}} & \cdots & 0 & 0 \\ \vdots & \vdots & \cdots & \vdots & \vdots \\ 0 & 0 & \ddots & \psi_{m-1,i^{(m-1)}} & -\psi m_{m,i^{(m)}} \\ 1 & 1 & \cdots & 1 & 1 \end{pmatrix} \tag{72}$$

Taking the derivations of (43) and $\sum_{(j,s)\in\Omega} \hat{\alpha}_{j,s} = 1$ with respect to $\hat{\alpha}_{\langle i \rangle, d_{\langle i \rangle}}, i \in \{2,\ldots,k\}$, we can get that $\hat{\Psi} \cdot \hat{\Phi}_i = \hat{Y}_i, i \in \{2,\ldots,k\}$, where $\hat{\Phi}_i \in \mathbb{R}^{m\times 1}$ and $\hat{Y}_i \in \mathbb{R}^{m\times 1}$,

$$\hat{\Phi}_i \triangleq \left( \frac{\partial \tilde{\alpha}_{\langle 1 \rangle,1}(y)}{\partial y}, \ldots, \frac{\partial \tilde{\alpha}_{\langle 1 \rangle,m}(y)}{\partial y} \right)' \bigg|_{y=\tilde{\alpha}_{\langle i \rangle, d_i}}, i \in \{2,\ldots,k\}, \tag{73}$$

$$\hat{\gamma}_i^{(l)} \triangleq \begin{cases} \gamma_{l+1,i} - \gamma_{l,i} &, i^{(l)} = i^{(l+1)} = i \\ -\gamma_{l,i} &, i^{(l)} = i, i^{(l+1)} \neq i \\ \gamma_{l+1,i} &, i^{(l)} \neq i, i^{(l+1)} = i \\ 0 &, i^{(l)} \neq i, i^{(l+1)} \neq i \end{cases}, l = 1,\ldots,m-1, i = 2,\ldots,k, \tag{74}$$

$$\hat{Y}_i \triangleq \left( \hat{\gamma}_1^{(l)}, \ldots, \hat{\gamma}_{m-1}^{(l)}, -1 \right)' \bigg|_{y=\tilde{\alpha}_{\langle i \rangle, d_i}}, i \in \{2,\ldots,k\}, \tag{75}$$

Where $'$ denotes the transpose operation of the matrix. In addition, by $\hat{\alpha}_{\langle i \rangle, d_i}$ and $\hat{\alpha}_{\langle i' \rangle, d_{i'}}$, $i \neq i', i, i' \in \{2,\ldots,k\}$ are independent, we have

$$\left. \frac{\partial M_{1,l,i,d_i}(\hat{\alpha}_{\langle 1 \rangle,l}, y)}{\partial y} \right|_{y=\hat{\alpha}_{\langle i \rangle, d_i}} = 0, i \in \{2,\ldots,k\}, i' \neq i. \tag{76}$$

And, following the allocation policy (24), (36)-(37), we can get that for $i \in \{2,\ldots,k\}$

$$\left.\frac{\partial M_{1,l,i^{(l)},d_{i^{(l)}}}(\hat{\alpha}_{\langle 1\rangle,l},y)}{\partial y}\right|_{y=\langle i\rangle^{(l)},d_{\langle i\rangle^{(l)}}} + \left.\frac{\partial M_{1,l,i^{(l)},d_{i^{(l)}}}(y,\hat{\alpha}_{\langle i\rangle^{(l)},d_{\langle i\rangle^{(l)}}})}{\partial y}\right|_{y=\langle 1\rangle,l} \left.\frac{\partial \hat{\alpha}_{\langle 1\rangle,l}(y)}{\partial y}\right|_{y=\hat{\alpha}_{\langle i\rangle^{(l)},d_{\langle i\rangle^{(l)}}}} = 0, l \in \{1,...,m\}.$$

in which if $i^{(l)} \neq \langle i \rangle$, $\left.\dfrac{\partial M_{1,l,i(l),d_{i(l)}}(\hat{\alpha}_{\langle 1\rangle,l},y)}{\partial y}\right|_{y=\hat{\alpha}_{\langle i(l)\rangle,d_{\langle i(l)\rangle}}} = 0$. Otherwise, $\exists s \in \{1,...,m\}$ such that the equality above does not hold, i.e.,

$$\left.\frac{\partial M_{1,s,i^{(s)},d_{i^{(s)}}}(\hat{\alpha}_{\langle 1\rangle,s},y)}{\partial y}\right|_{y=\langle i\rangle,d_{\langle i\rangle}} + \left.\frac{\partial M_{1,s,i^{(s)},d_{i^{(s)}}}(y,\hat{\alpha}_{\langle i\rangle^{(s)},d_{\langle i\rangle^{(s)}}})}{\partial y}\right|_{y=\langle 1\rangle,s} \left.\frac{\partial \hat{\alpha}_{\langle 1\rangle,s}(y)}{\partial y}\right|_{y=\hat{\alpha}_{\langle i\rangle,d_{\langle i\rangle}}} = 0, l \in \{1,...,m\}.$$

Following the RAODA, (36)-(37), alternative-scenario pair $(\langle i \rangle, d_{\langle i \rangle})$ will receiving simulation replications and alternative-scenario pair $(\langle 1 \rangle, s)$ will stop being allocated simulation budget before the inequality sign above is no longer hold, which contradicts with $\left(\alpha_{1,1}^{(t)}, ..., \alpha_{1,m}^{(t)} \alpha_{2,d_2}^{(t)}, ..., \alpha_{k,d_k}^{(t)}\right)$ converges to $(\hat{\alpha}_{1,1}, ..., \hat{\alpha}_{1,m}, \hat{\alpha}_{2,d_2}, ..., \hat{\alpha}_{k,d_k})$. Thus, we have $\widehat{K} \cdot \widehat{\Phi}_i = \widehat{G}_i, i \in \{2,...,k\}$, where $\widehat{K} \in \mathbb{R}^{(k-1)\times(k-1)}$, $\widehat{G}_i \in \mathbb{R}^{(k-1)\times 1}$, and the diagonal matrix

$$\widehat{K} \triangleq \text{diag}\left(\psi_{1,i^{(1)}}, \psi_{2,i^{(2)}}, ..., \psi_{m,i^{(m)}}\right), \tag{77}$$

and

$$\hat{g}_l^{(i)} \triangleq \begin{cases} -\gamma_{l,i+1} & l^{(i+1)} = \langle l \rangle \\ 0 & l^{(i+1)} \neq \langle l \rangle \end{cases}, l=1,2,\cdots,m, i=1,2,...,k-1, \tag{78}$$

$$G_l \triangleq \left(g_1^{(i)}, g_2^{(i)}, ..., g_m^{(i)}\right)', \tag{79}$$

in which $'$ denotes the transpose operation of the matrix. According to the Lemma 3 in Zhang et al., (2022). we have $\psi_{l,i} > 0, l=1,...,m, i=2,...,k$, which indicates $\det(\widehat{K}) > 0$ and $\widetilde{K}$ is invertible,

$$K^{-1} \triangleq \text{diag}\left(\psi_{1,i^{(1)}}^{-1}, \psi_{2,i^{(2)}}^{-1}, ..., \psi_{m,i^{(m)}}^{-1}\right), \tag{80}$$

So, we have $\Upsilon_l = \widehat{\Psi} \cdot \widehat{K}^{-1} \widehat{G}_i, i=2,...,k$, which leads to

$$\sum_{\{l|(\langle i\rangle,d_{i^{(l)}})=(\langle i\rangle,d_{\langle i\rangle})\}} \frac{\partial M_{1,l,i^{(l)},d_{i^{(l)}}}(\hat{\alpha}_{\langle 1\rangle,l},y)/\partial y|_{y=\hat{\alpha}_{i^{(l)},d_{i^{(l)}}}}}{\partial M_{1,l^{(i)},i,d_i}(y,\hat{\alpha}_{i^{(l)},d_{i^{(l)}}})/\partial y|_{y=\hat{\alpha}_{\langle 1\rangle,l}}} = 1, i=2,...,k. \tag{81}$$

Thus，we can get that

$$\frac{\left(\alpha^*_{\langle i\rangle,d_{\langle i\rangle}}\right)^2}{\sigma^2_{\langle i\rangle,d_{\langle i\rangle}}} = \sum_{\{l|(\langle i\rangle,d_{\langle i\rangle^l})=(\langle i\rangle,d_{\langle i\rangle})\}} \frac{\left(\alpha^*_{\langle 1\rangle,l}\right)^2}{\sigma^2_{\langle 1\rangle,l}}, l=1,\ldots,m. \tag{82}$$

Since $\bigcup_{l=2,\ldots,k}\{l|(\langle i\rangle^l, d_{\langle i\rangle^l}) = (\langle i\rangle, d_{\langle i\rangle})\} = \{1,\ldots,m\}$, $\bigcap_{l=2,\ldots,k}\{l|(\langle i\rangle^l, d_{\langle i\rangle^l}) = (\langle i\rangle, d_{\langle i\rangle})\} = \emptyset$, we have (40) must hold. ∎

If (38)-(40) has the unique solution, Theorem 2 can be adapted to prove that the allocation (23) possesses the following asymptotic optimal property.

**Theorem 3.** Suppose (38)-(40) has a unique solution, then

$$\lim_{t\to\infty} \alpha^{(t)}_{j,s} \to \alpha^*_{j,s} \quad a.s. \tag{83}$$

where $\alpha^{(t)}_{j,s} = t_{j,s}/t$, $\forall (j,s) \in \Omega$.

*Proof:* We assume that $\left(\alpha^{(t)}_{1,1},\ldots,\alpha^{(t)}_{1,m}\alpha^{(t)}_{2,d_2},\ldots,\alpha^{(t)}_{k,d_k}\right)$ not converges to $(\alpha^*_{1,1},\ldots,\alpha^*_{1,m}\alpha^*_{2,d_2},\ldots,\alpha^*_{k,d_k})$, by Bolzano-Weierstrass theorem, there exist a subsequence of $\left(\alpha^{(t)}_{1,1},\ldots,\alpha^{(t)}_{1,m}\alpha^{(t)}_{2,d_2},\ldots,\alpha^{(t)}_{k,d_k}\right)$ converges to $(\hat{\alpha}_{1,1},\ldots,\hat{\alpha}_{1,m},\hat{\alpha}_{2,d_2},\ldots,\hat{\alpha}_{k,d_k}) \neq (\alpha^*_{1,1},\ldots,\alpha^*_{1,m}\alpha^*_{2,d_2},\ldots,\alpha^*_{k,d_k})$, Without loss of generality, we can suppose that $\left(\alpha^{(t)}_{1,1},\ldots,\alpha^{(t)}_{1,m}\alpha^{(t)}_{2,d_2},\ldots,\alpha^{(t)}_{k,d_k}\right)$ converges to $(\hat{\alpha}_{1,1},\ldots,\hat{\alpha}_{1,m},\hat{\alpha}_{2,d_2},\ldots,\hat{\alpha}_{k,d_k})$. Following the Theorem 2 and the premise that (38)-(40) has a unique solution, we have $(\hat{\alpha}_{1,1},\ldots,\hat{\alpha}_{1,m},\hat{\alpha}_{2,d_2},\ldots,\hat{\alpha}_{k,d_k}) = (\alpha^*_{1,1},\ldots,\alpha^*_{1,m}\alpha^*_{2,d_2},\ldots,\alpha^*_{k,d_k})$, which contradicts the assumption that $\left(\alpha^{(t)}_{1,1},\ldots,\alpha^{(t)}_{1,m}\alpha^{(t)}_{2,d_2},\ldots,\alpha^{(t)}_{k,d_k}\right)$ converges to $(\hat{\alpha}_{1,1},\ldots,\hat{\alpha}_{1,m},\hat{\alpha}_{2,d_2},\ldots,\hat{\alpha}_{k,d_k})$. Therefore, $\left(\alpha^{(t)}_{1,1},\ldots,\alpha^{(t)}_{1,m}\alpha^{(t)}_{2,d_2},\ldots,\alpha^{(t)}_{k,d_k}\right)$ converges to $(\alpha^*_{1,1},\ldots,\alpha^*_{1,m}\alpha^*_{2,d_2},\ldots,\alpha^*_{k,d_k})$. ∎

## 5. Numerical Experiments

In this section, we present the RAODA procedure for selecting the best alternative problem with input uncertainty based on RAODA (24) and test the efficiency of the proposed the RAODA procedure by a series of numerical experiments.

---
**RAODA Procedure**

**Input:** Total number of alternatives $k$, the number of scenarios for each alternative $m$, total number of simulation replications $T$, and the initial simulation replication number $n_0$.

> **Initialize:** $t = 0$ and perform $n_0$ simulation replications for each design-scenario pair.
> **while** $t < T - n_0 * k * m$, **Do:**
>     **Update:** Calculate the posterior variances $\left(\sigma_{i,j}^{(t)}\right)^2$ and posterior means $\mu_{i,j}^{(t)}$ with sample means $\bar{X}_{i,d}^{(t)}$ and sample variances $\left(\bar{\sigma}_{i,d}^{(t)}\right)^2$, according to equation (16).
>     **Allocation:** Determine the allocated alternative-scenario pair $\hat{a}_{t+1}(\mathcal{F}_t)$ according to equation (22)-(24).
>     **Simulation:** Run an additional simulation for the allocated alternative.
>     $t \leftarrow t + 1$;
> **end while.**

The RAODA procedure is compared with the following allocation rules:

●Robust optimal computing budget allocation procedure (ROCBA) for the selection of the best with input uncertainty (Gao et al., 2017).

● Equal allocation procedure (EA) is the simplest allocation algorithm, and it is often used as a benchmark for comparison with other procedures. The procedure equally allocates simulation budget to estimate the performance of each alternative, i.e., roughly $T/t$ samples for each alternative-scenario pair.

● Proportional to variance procedure (PTV) takes $n_0$ independent samples of each alternative during the initial stage and uses the newly updated sample variances to guide the allocation of additional increments to the simulation budget for each period, i.e., $\alpha_{i,d}/\alpha_{j,s} = s_{i,d}^2/s_{j,s}^2$, $\forall (i,d), (j,s) \in \{1,2,\ldots,k\} \times \{1,2,\ldots,m\}$ ( Rinott 1978). Notably, the PTV allocation rule is equivalent to the EA allocation rule when the variances of each alternative are equal.

In all numerical examples, $PCS_W$ is applied as an objective function of the sampling budget $T = 6000$ in each experiment which is used to measure the statistical efficiency of the sampling process estimated by Times $= 10000$ independent experiments. The initial number of simulation replications $n_0 = 40$, and incremental replications is also $\Delta = 1$.

To compare RAODA, ROCBA, EA, PTV, we provide the following the selection of the best with input uncertainty:

**Experiment 1:** We have $k = 10$ alternatives, and each alternative has $m = 5$ scenarios. the prior distribution of the performance $\mu_{i,d}$ of each alternative is Gaussian distribution with prior mean $\mu_{i,d}^{(0)} = 0$ and prior variance $\left(\sigma_{1,d}^{(0)}\right)^2 = 0.01$ , $\left(\sigma_{j,d}^{(0)}\right)^2 = 0.02$, $j \neq 1, i, j \in \{1,2,\ldots,k\}, d \in \{1,2,\ldots,m\}$ and the simulation replications are drawn independently from a

Gaussian distribution $N\left(\mu_{i,d},(\sigma_{i,d})^2\right)$, where $(\sigma_{i,d})^2 = 1$ for $(i,d) \in \{1,2,\ldots,k\} \times \{1,2,\ldots,m\}$.

**Experiment 2:** We have $k = 10$ alternatives, and each alternative has $m = 5$ scenarios. the prior distribution of the performance $\mu_{i,d}$ of each alternative is Gaussian distribution with prior mean $\mu_{i,d}^{(0)} = 0$ and prior variance $\left(\sigma_{i,d}^{(0)}\right)^2 = (3 - 0.1 * (i + d))^2$ and the simulation replications are drawn independently from a Gaussian distribution $N\left(\mu_{i,d},(\sigma_{i,d})^2\right)$, where $(\sigma_{i,d})^2 = 8^2$ for $(i,d) \in \{1,2,\ldots,k\} \times \{1,2,\ldots,m\}$.

**Experiment 3:** We have $k = 10$ alternatives, and each alternative has $m = 5$ scenarios. the prior distribution of the performance $\mu_{i,d}$ of each alternative is Gaussian distribution with prior mean $\mu_{i,d}^{(0)} = 0$ and prior variance $\left(\sigma_{i,d}^{(0)}\right)^2 = (i + d)^2$ and the simulation replications are drawn independently from a Gaussian distribution $N\left(\mu_{i,d},(\sigma_{i,d})^2\right)$, where $(\sigma_{i,d})^2 = 8^2$ for $(i,d) \in \{1,2,\ldots,k\} \times \{1,2,\ldots,m\}$.

From Figures 2-4, the RAODA procedure proposed in this paper has the best performance among all numerical experiments, the ROCBA procedure has the second-best performance, while the PTV and EA procedures are not as good as the ROCBA procedure. In Experiments 1-3, to achieve the PCS at 0.45, 0.65, and 0.93, the ROCBA procedure consumes an extra simulation budget of 856, 510, and 2724 respectively, compared with the RAODA procedure. The analysis of the above experimental results shows that the proposed RAODA procedure significantly improves the simulation efficiency of the optimal selection problem under input uncertainty.

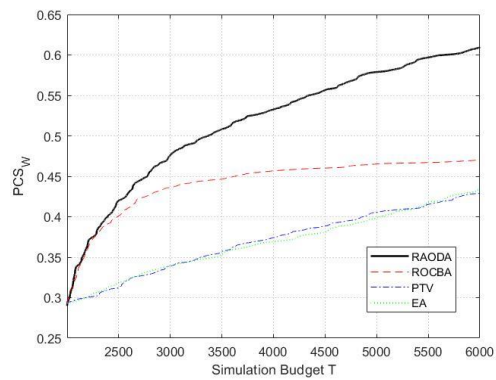

Fig. 2. PCS of the four procedures for experiment 1

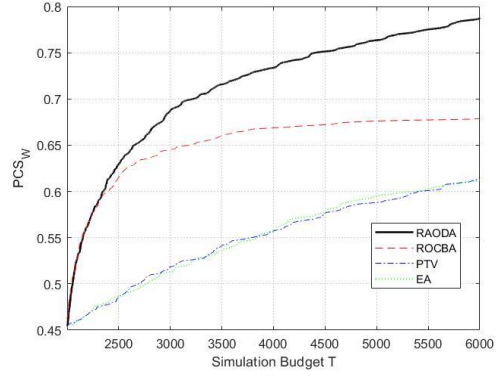

Fig. 3. PCS of the four procedures for experiment 2

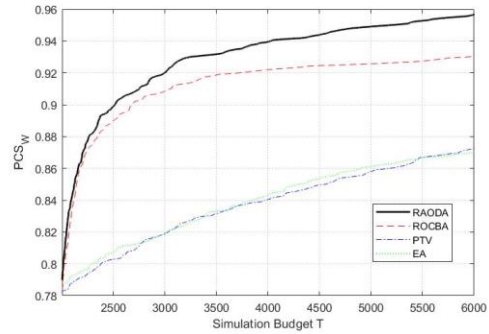

Fig. 4. PCS of the four procedures for experiment 3

**6. Conclusions**

This paper considers the selection of the best with input uncertainty problem. In a Bayesian framework, we describe the simulation budget allocation decision process as a stochastic control problem based on a stochastic dynamic programming approach and equivalently transform it into an MDP. Following the value function approximation, we derive the asymptotically optimal dynamic allocation policy (RAODA) with the objective of maximizing the posterior probability of the correct selection of the best ($PCS_W$) and show that the RAODA can asymptotically achieve consistency and optimality. In addition, numerical experiments on three cases show that the RAODA procedure significantly improves the simulation efficiency of the robust ranking and selection problem compared with existing algorithms.